\documentclass{elsart1p}
\usepackage[french]{babel}
\usepackage[T1]{fontenc}
\usepackage{amsmath}
\usepackage{amssymb}
\usepackage[dvips]{graphicx}
\usepackage{float}
\usepackage{bbold}
\usepackage{cancel}
\usepackage{dsfont}
\usepackage{amsfonts}

\newtheorem{dfn}[thm]{\textsf{Définition}}

\newcommand{\N}{\mathds{N}}
\newcommand{\R}{\mathds{R}}

\begin{document}
\begin{frontmatter}

\title{Solutions globales pour des équations de Schrödinger sur-critiques en toutes dimensions}
\author[South]{A.Poiret},
\ead{aurelien.poiret@math.u-psud.fr}
\address[South]{Faculté des sciences d'Orsay, Département de mathématiques, Bâtiment 430 Bureau 106, France}

\begin{abstract} 
Dans \cite{moi}, on a expliqué comment construire un grand nombre de solutions globales pour l'équation de Schrödinger cubique en dimension 3 avec des données initiales dans $ L^2( \mathds{R}^3) $. Les arguments de bases étant vraie en dimensions plus grandes que 2, nous pouvons adapter la preuve dans ces cas là. On explique dans cet article comment utiliser l'effet régularisant pour prouver un théorème analogue en toutes dimensions, en particulier en dimension 1. Le gain de régularité est plus faible mais l'on peut choisir une base de fonctions propres quelconques et des variables aléatoires autres que gaussiennes.
\end{abstract}


\end{frontmatter}

\newpage
Dans cet article, on considère les équations de Schrödinger suivantes :
\begin{equation} \label{2schrodinger} \tag{$NLS$}
  \left\{
      \begin{aligned}
       & i \frac{ \partial \tilde{ u }  }{ \partial t } +  \Delta  \tilde{ u }   =   K | \tilde{u} |^{p-1} \tilde{ u },
       \\ &  \tilde{u}(0,x)=u_0(x),
      \end{aligned}
    \right.
\end{equation}
où $ K \in \lbrace -1,1 \rbrace $ et p désigne un entier impair.
\\\\Dans \cite{moi}, on donne une méthode pour construire des solutions globales pour les équations de Schrödinger dont le nombre de dérivée sur-critiques est inférieure à $ \frac{1}{2} $, en dimension plus grande que 2. L'idée de la preuve est de rendre la donnée initiale aléatoire et d'utiliser des estimées bilinéaires de type Bourgain. Ici, on propose de compléter ce résultat, en particulier en établissant le théorème en dimension 1.
\\\\En dimension 1, dans \cite{burq6}, il est prouvé que l'effet régularisant permet de gagner $ \frac{1}{2} - \frac{2}{p-1} $ dérivée sur la donné initiale. Cela signifie que des estimées linéaires sont suffisantes pour établir le résultat en dimension 1. Cette méthode est très spécifique à la dimension 1 et ne peut être généralisée directement en dimension plus grande. Néanmoins, pour $ p \geq 5 $ et $ u_0 \in \overline{H}^{(d-1)/2} (  \mathds{R}^d) $, le même schéma de preuve permet de gagner le nombre de dérivée manquant.
\\\\Ce résultat est très intéressant car il n'est plus nécessaire de supposer que les fonctions propres soient les fonctions tenseurs. Une base de fonctions propres quelconques est satisfaisante et le théorème est vérifié pour un plus grand nombre de mesure de probabilité.
\\De plus, on propose une preuve du théorème dans un cadre plus général que des variables aléatoires gaussiennes ou Bernoulli.
\section{Introduction et notations}
En dimension d'espace d quelconque, on pose $ H = -\Delta + x^2 $ l'oscillateur harmonique. On note $ \lambda_n ^2 $ les valeurs propres et  $  h_n $ les fonctions propres de H que l'on indexe par $ n \in \mathds{N} $. On a donc
\begin{equation*}
H h_n = \lambda_n^2 h_n , \ \forall n \in \N.
\end{equation*}
On note $ W^{s,p}(\R^d) $ et $ H^s(\R^d) $ les espaces de Sobolev usuels. Puis, on définit les espaces de Sobolev harmoniques.
\begin{dfn} L'espace $ \overline{H}^s( \R^d ) $ est défini comme la fermeture de l'espace de Schwartz pour la norme
\begin{equation*}
|| u || _{  \overline{H} ^s(  \R^d )   } = || H^{s/2} u ||_ {L^2(   \R^d ) }.
\end{equation*}
\end{dfn}
\begin{dfn}
De manière similaire, l'espace $ \overline{W}^{s,p}( \R^d ) $ est défini comme la fermeture de l'espace de Schwartz pour la norme
\begin{equation*}
|| u || _{  \overline{W} ^{s,p}(  \R^d )   } = || H^{s/2} u ||_ {L^p(   \R^d ) }.
\end{equation*}
\end{dfn}
Dans \cite{sobolev}, nous pouvons trouver la proposition suivante :
\begin{prop} \label{comparaison}
Pour tous $ 1 < p < \infty $, $ s \geq 0 $, il existe une constante $ C > 0 $ telle que 
\begin{equation*} 
\frac{1}{C} || u || _{  \overline{W}^{s,p}( \R^d ) } \leq  || \nabla^s u ||_{L^p(\R^d)} + || <x>^s u ||_{L^p(\R^d ) }  \leq C || u || _{  \overline{W}^{s,p}( \R ^d ) }.
\end{equation*}
\end{prop}
Ensuite, soit $ (\Omega , A , P)  $ un espace de probabilité, $ (g_n ( \omega ) )_{n \in \mathds{N}} $ une suite de variable aléatoires indépendantes et définissons les conditions suivantes :
\begin{equation}  \label{2hypothese1} \tag{$H_\gamma $}
\boxed{ \exists \ \gamma , C, c > 0 \ / \ \forall n \in \mathds{N} \mbox{ and } \rho  \in \R, \   \int_\rho^\infty dP_{g_n} + \int_{- \infty }^{-\rho} dP_{g_n}   \leq C e ^{-c |\rho|^\gamma}}
\end{equation}
\begin{equation} \label{2hypothese2} \tag{$H_{E_1}$} 
 \boxed{ \forall p \in \N \mbox{ et } n \in \mathds{N}, \  E( g_n ^{2p+1} )   = 0}   
\end{equation}
\begin{equation}   \label{2hypothese2bis} \tag{$H_{E_2}$} 
 \boxed{ \forall n \in \mathds{N}, \  E( g_n )   = 0} 
\end{equation}
\begin{equation}  \label{2hypothese3} \tag{$H_{01}$}
\boxed{ \forall \rho > 0 \mbox{ et } n \in \mathds{N} , \  P ( |g_n| < \rho  )   > 0 } 
\end{equation}
\begin{equation}  \label{2hypothese4} \tag{$H_{02}$}
\boxed{ \exists c > 0 / \ \forall n \in \mathds{N} , \  E(|g_n|^2)   \geq c }  
\end{equation}
Nous avons facilement le lemme suivant:
\begin{lem} \label{2moments}
Sous l'hypothèse (\ref{2hypothese1}), il existe une constante $ C> 0 $ telle que pour tout $ n \in \mathds{N} $,
\begin{equation*}
 E(|g_n|^2)^2 \leq E(|g_n|^4) \leq C.
\end{equation*}
\end{lem}
\textit{Preuve.}  $ E(|g_n|^4) = 4 \int_0^\infty \rho^3 P \left( \omega \in \Omega / |g_n(\omega)| \geq \rho \right) d \rho  \leq 4C \int_0^\infty \rho^3  e^{-c \rho^\gamma} d \rho  < \infty. $ \hfill $ \boxtimes $
\\\\Soit $ u_0 \in \overline{H}^{  \frac{d-1}{2} } (  \R ^d ) $, c'est à dire
\begin{equation*}
u_0(x) = \sum_{n \in \mathds{N} } c_\lambda h_n(x)     \ \mbox{avec} \  \sum_{n \in \mathds{N} }  \lambda_n ^{d-1} |c_n|^2 < \infty.
\end{equation*}
Considérons l'application $ \omega \longrightarrow u_0^\omega  $ de $ ( \Omega , A , P ) $ dans $ \overline{H}^{  \frac{d-1}{2} } (  \mathds{R} ^d )  $  que l'on équipe de sa tribu borélienne, définie par $ u_0^\omega  = \displaystyle{ \sum_{n \in \mathds{N}} } c_n g_n( \omega ) h_n(x) $.
\\\\Grâce au lemme \ref{2moments}, nous pouvons facilement vérifier que l'application $ \omega \longrightarrow u_0^\omega $ est dans $ L^2 ( \Omega , \overline{H}^{  \frac{d-1}{2} } (  \mathds{R} ^d ) )  $. Enfin, on définit $ \mu $  comme la loi de la variable aléatoire $ \omega \rightarrow u_0(\omega,.) $ et nous pouvons donc appliquer le théorème de transfert suivant :
\begin{equation}
P( \omega \in \Omega / \Psi ( u_0^\omega ) \in A ) = \mu ( u_0 \in \overline{H}^{  \frac{d-1}{2} } (  \R ^d ) / \Psi (u_0) \in A   ),
\end{equation}
pour toute fonction mesurable $ \Psi $ et $ A $ ensemble mesurable.
\\\\Pour pouvoir énoncer les théorèmes de ce papier, on introduit les deux définitions suivantes :
\begin{dfn}
Soit $ (q , r ) \in [2, \infty ] $ alors $ ( q , r ) $ est dit admissible si
\begin{align*}
 (d,q,r) \neq (2,2,\infty)  \mbox{ et }   \frac{2}{q} = \frac{d}{2} - \frac{d}{r}.
\end{align*}
\end{dfn}
\begin{dfn} Pour $ s \in \R $ et $ T > 0 $, on définit
\begin{align*}
 &  X^s = \underset{ (q,r) \ admissible }{ \bigcap } L^q( \mathds{R} , W^{s,r}(\mathds{R}^d)),
\\ & X_T^s = \underset{ (q,r) \ admissible }{ \bigcap } L^q( [-T,T] , W^{s,r}(\mathds{R}^d)).
\end{align*}
\end{dfn}
Dans cet article, sous les hypothèses (\ref{2hypothese1}), (\ref{2hypothese2}), (\ref{2hypothese3}) et (\ref{2hypothese4}) \\  ou (\ref{2hypothese1}), (\ref{2hypothese2bis}), (\ref{2hypothese3}) et (\ref{2hypothese4}), on propose de démontrer les théorèmes suivants :
\begin{thm} \label{2thm1}
Soit $ u_0 \in \overline{H}^{  \frac{d-1}{2} } (  \mathds{R} ^d ) $ alors il existe $ s \in ]  \frac{d}{2} - \frac{2}{p-1}  , \frac{d}{2} [ $ et un ensemble $ \Omega ' \subset \Omega  $ tels que les conditions suivantes soient réalisées : 
\\i) $ P (  \Omega '  ) > 0 $.
\\ii) Pour tout élément $ \omega \in \Omega'$, il existe une unique solution globale $ \tilde{u} $ à l'équation (\ref{2schrodinger}) dans l'espace $ e^{it\Delta} u_0 (\omega,.) + X ^s  $ avec donnée initiale $ u_0 (\omega,.) $.
\\iii) Pour tout élément $ \omega \in \Omega ' $, il existe $ L^+ \in \overline{H}^s(\mathds{R}^d) $ et $ L_- \in \overline{H}^s(\mathds{R}^d) $ telles que 
\begin{align*}
&  \lim _{ t \rightarrow \infty  } || \tilde{u}(t) - e^{it\Delta} u_0 (\omega,.) - e^{it\Delta} L^+ ||_{  H^s(\mathds{R}^d) } = 0,
\\ & \lim _{ t \rightarrow - \infty  } ||  \tilde{u}(t) - e^{it\Delta} u_0 (\omega,.) - e^{it\Delta} L_-  ||_{   H^s(\mathds{R}^d) } = 0.
\end{align*}
De plus, si $ u_0 \notin \overline{H}^s(\mathds{R}^d) $ alors $ P \left( \omega \in \Omega / u_0 (\omega,.) \in H^s( \mathds{R}^d )  \right) = 0 $.
\end{thm}
\begin{thm} \label{2thm3}
Soit $ u_0 \in \overline{H}^{  \frac{d-1}{2} } (  \mathds{R} ^d ) $ alors il existe $ s \in ]  \frac{d}{2} - \frac{2}{p-1} , \frac{d}{2} [ $  tel que pour tout $ \omega \in \Omega $, il existe $ T_\omega $ et une unique solution à l'équation (\ref{2schrodinger}) dans l'espace $ e^{it\Delta} u_0 (\omega,.) + X_{T_\omega} ^s $ avec donnée initiale $ u_0 (\omega,.)  $.
\\\\Plus précisément, il existe $ C,c, \delta > 0 $ et pour tout temps $ 0 < T < \infty $, un ensemble $ \Omega_T $ tels que
\begin{equation*}
P(  \Omega _T  ) \geq 1-C e^{  -c/  \arctan ( 2 T )^\delta },
\end{equation*} 
et tel que pour tout élément $ \omega \in \Omega_T $, il existe une unique solution à l'équation (\ref{2schrodinger}) avec donnée initiale $ u_0 (\omega,.)  $ dans un espace continûment inclus dans $ C ^0 ( [ - T , T ]  ,  H^{  \frac{d-1}{2} } (  \mathds{R} ^d ) )$.
\end{thm}
\begin{thm} \label{2thm2} Si de plus, pour tout $ n \in \mathds{N} $, $ g_n $ a une distribution symétrique, alors  
\begin{align*}
 \lim_{ \eta \rightarrow 0  } \  \mu  \left( u_0 \in \overline{H}^{  \frac{d-1}{2} } (  \mathds{R} ^d )   / \mbox{ on ait existence globale et scattering } | \ || u_0 ||_{\overline{H}^{  \frac{d-1}{2} } (  \mathds{R} ^d )} \leq \eta  \right) = 1.
\end{align*}
\end{thm}
\section{Quelques rappels préliminaires}
Dans cette section, on rappelle les estimées de Strichartz pour l'oscillateur harmonique ainsi que la propriété fondamentale de la transformation de lentille. Les preuves peuvent être trouvées dans \cite{moi}.
\subsection{Estimées de Strichartz pour l'oscillateur harmonique}
\begin{dfn} Pour $ s \in \R $ et $ T \geq 0 $, on définit
\begin{align*}
\overline{X }_T^s = \underset{ (q,r) \ admissible }{ \bigcap } L^q( [-T,T] , \overline{W}^{s,r}(\mathds{R}^d)).
\end{align*}
\end{dfn}
Alors, nous avons les propositions suivantes :
\begin{prop} \label{2Stricharz1} Pour tout $ T \geq 0 $, il existe une constante $ C_T > 0 $ telle que pour tout $ u \in \overline{H}^s(\mathds{R}^d) $,
\begin{equation*}
 || e^{-itH} u  ||_{\overline{X }_T^s} \leq C_T || u ||_{  \overline{H}^s(\mathds{R}^d) }.
\end{equation*}
\end{prop}
\begin{prop} \label{2Stricharz2} Pour tout $ T \geq 0 $, il existe une constante $ C_T > 0 $ telle que pour tout $ ( q, r ) $ admissible et $ F \in L^{q'}( [-T,T], \overline{W}^{s,r'} (\mathds{R}^d)) $,
\begin{equation*}
 \bigg| \bigg| \int _0^t e^{-i(t-s)H} F(s) ds   \bigg| \bigg| _{  \overline{X}^s_T} \leq C_T || F ||_{  L^{q'}( [-T,T], \overline{W}^{s,r'} (\mathds{R}^d)) }.
\end{equation*}
\end{prop}
\subsection{La transformation de lentille}
\begin{dfn}
Pour $ u(t,x)$ une fonction mesurable avec $ t \in \R $ et $ x \in \mathds{R}^d $, on définit $ \tilde{u}(t,x) $ de la façon suivante :
\begin{equation*}
\tilde{u}(t,x) =  \left( \frac{1}{\sqrt{1+4t^2}}  \right) ^{d/2} \times u \left( \frac{1}{2} \arctan(2t) , \frac{x}{\sqrt{1+4t^2} } \right)   \times e^{ \frac{ix^2t}{1+4t^2}  }.
\end{equation*}
\end{dfn}
En analogie à la proposition 22 de \cite{moi}, on peut obtenir le résultat suivant :
\begin{prop} \label{2comparaisonbis} Soit $ s \geq 0$ alors il existe une constante $ C > 0 $ telle que pour tout $ p \in [1,+\infty ] $ et $ q \in [1,+\infty ]  $ vérifiant $ \frac{2}{p} + \frac{d}{q} - \frac{d}{2} \leq 0   $, pour tout $ T \in ] 0 , \infty ] $ et $ u \in L^p( [  - \frac{1}{2} \arctan(2T)   , \frac{1}{2} \arctan(2T) ] , \overline{W}^{s,q} ( \mathds{R}^d ) ) $, on a 
\begin{align*}
|| \tilde{u} ||_{ L^p( [-T,T] , W^{s,q} ( \mathds{R}^d )  )  } \leq C || u || _{L^p( [  - \frac{1}{2} \arctan(2T)   , \frac{1}{2} \arctan(2T) ] , \overline{W}^{s,q} ( \mathds{R}^d ) )   }.
\end{align*}
\end{prop}
\section{L'effet régularisant pour l'oscillateur harmonique}
On commence par donner une preuve de l'effet régularisant de \cite{yajima1} et \cite{yajima2} en utilisant une méthode de Doï. Cet effet régularisant se révèlera fondamental pour appliquer le théorème de point fixe de Picard. L'objectif de cette partie est donc de prouver le théorème suivant :
\begin{thm}  
Soit $ \epsilon \in ] 0 , \frac{1}{2} [ $ alors il existe une constante $ C > 0 $ telle que pour tout $ u_0 \in L^2(\R^d) $,
\begin{equation} \label{2effectregularisant1} \bigg| \bigg| \frac{1}{<x>^{1/2-\epsilon} }\sqrt H ^{1/2-2\epsilon} e^{itH} u_0 \bigg| \bigg|_{L^2([-2 \pi,2\pi]*\R^d)} \leq C ||u_0||_{L^2(\R^d)},
\end{equation}
et pour tout $ u_0 \in  \overline{H}^{   \frac{d-1}{2} } (\R^d) $,
\begin{equation} \label{2effectregularisant2} \bigg| \bigg| \frac{1}{<x>^{1/2-\epsilon} } |\nabla  | ^{d/2-2\epsilon} e^{itH} u_0 \bigg| \bigg|_{L^2([-2 \pi,2\pi]*\R^d)} \leq C ||u_0||_{\overline{H}^{   \frac{d-1}{2} } (\R^d)}.
\end{equation}
\end{thm}
\subsection{Quelques résultats préliminaires}
On commence par établir 4 lemmes préliminaires.
\begin{lem}  \label{2inter}
Soit $ a \in C^\infty(\R^d, \R^d) \cap L^\infty( \R^d, \R^d) $ telle que $ \nabla a  \in L^\infty (\R^d, M_d(\R^d)) $ alors il existe une constante $ C> 0 $ telle que pour tout $ u \in H^{1/2}(\R^d) $,
\begin{equation*}
 \bigg| \int_{\R^d} a(x) . \nabla u (x)   \overline{u} (x)  \ dx \bigg|  \leq C ||u||^2_{H^{1/2}(\R^d)}.
\end{equation*}
\end{lem}
 \textit{Preuve.} On définit
\begin{equation*}
b(u,v ) = \int_{\R^d} a(x) . \nabla u(x) v(x) \ dx.
\end{equation*}
Alors, clairement, on a  
\begin{equation*}
| b(u,v ) | \leq C ||u||_{H^1(\R^d)} || v ||_{L^2 (\R^d)},
\end{equation*}
ainsi que 
\begin{equation*}
 | b(u,v ) | = \bigg| \int_{\R^d } u (x) \sum_{i=1}^d  \partial_i \left( a_i (x) v (x) \right) \ dx \bigg|   \leq C  ||u||_{L^2(\R^d)} ||v||_{H^1(\R^d)}.
 \end{equation*}
Par conséquent, par interpolation, pour tout $ s \in [0,1]$, il existe une constante $ C>0$ telle que \begin{equation*} 
|b(u,v)| \leq C ||u||_{H^s(\R^d)} ||v||_{H^{1-s}(\R^d)}.
\end{equation*}
Le lemme est donc démontré en choisissant $ s = 1/2 $. \hfill $ \boxtimes $
\begin{lem}  \label{2interbis} Soit $ a \in C^\infty( \R^d, \R^d) $ telle que  $ |a(x)| \leq |x|^{2\epsilon}$ et $ |\nabla a (x) | \leq 1 $ alors il existe une constante $ C> 0 $ telle que pour tout $ u  \in \overline{H}^{1/2+\epsilon}(\R^d ) $,
\begin{equation*}
\bigg| \int_{\R^d} a(x). \nabla u( x) \overline{u}(x) \ dx   \bigg| \leq C ||u||_{\overline{H}^{1/2+\epsilon}(\R^d)}^2.
\end{equation*}
\end{lem}
\textit{Preuve.} Il s'agit essentiellement de la même preuve que le lemme \ref{2inter}.
\\\\On définit 
\begin{equation*}
b(u,v ) = \int_{\R^d} a(x) . \nabla u(x) v(x) \ dx.
\end{equation*}
Alors, on trouve
\begin{equation*} 
| b(u,v ) | \leq C ||u||_{H^1(\R^d)} || v ||_{\overline{H}^{2\epsilon}(\R^d) } \leq C ||u||_{  \overline{H}^1(\R^d)} || v ||_{\overline{H}^{2\epsilon}(\R^d) },
\end{equation*} 
ainsi que
\begin{align*}  | b(u,v ) | = \bigg| \int_{\R^d } u (x) \sum_{i=1}^d  \partial_i  \left( a_i (x) v (x)  \right) \ dx \bigg| & \leq C ||u||_{\overline{H}^{2\epsilon}(\R^d)} ||v||_{H^1(\R^d)}
\\ & \leq C ||u||_{\overline{H}^{2\epsilon}(\R^d)} ||v||_{ \overline{H}^1(\R^d)}.
\end{align*}
Puis, par interpolation, il existe une constante $ C> 0 $ telle que pour tout $ s \in [0,1]$,  
\begin{equation*} 
|b(u,v)| \leq C ||u||_{\overline{H}^{(1-2\epsilon )s+2\epsilon}(\R^d)} ||v||_{\overline{H}^{1-s(1-2\epsilon)}(\R^d)}.
\end{equation*} 
Le lemme est donc prouvé en prenant $ s = 1/2 $. \hfill $ \boxtimes $
\begin{lem}  \label{2commutateur}
Soient $ s_1 $ et $ s_2 $ deux réels.
\\ - Si $  \max( s_2,s_1+ s_2 ) \leq 1 $ alors il existe une constante $ C > 0 $ telle que pour tout $ u \in L^2(\R^d) $,
\begin{align*}
|| \ [\sqrt H^{s_1+s_2};<x>^{-s_1}]u ||_{L^2(\R^d) } \leq C ||u ||_ {L^2(\R^d)},
 \end{align*}
- Si $ s_2 \geq -1 $ alors il existe une constante $ C > 0 $ telle que pour tout $ u \in  H^{s_1-1} (\R^d) $,
\begin{align*}
 || \ [  \ |\nabla | ^{s_1};<x>^{-s_2}]u ||_{L^2(\R^d) } \leq C ||u ||_ {H^{s_1-1} (\R^d)},
\end{align*}
- Si $ s_2 \leq 1 $ alors il existe une constante $ C > 0 $ telle que pour tout $ u \in  \overline{H}^{s_1-s_2} (\R^d) $,
\begin{align*}
 || \ [  \sqrt{H}^{s_1};<x>^{-s_2}]u ||_{L^2(\R^d) } \leq C ||u ||_ {\overline{H}^{s_1-s_2} (\R^d)}.
\end{align*}
\end{lem}
\textit{Preuve.} Pour évaluer la régularité du commutateur, on utilise le calcul pseudo-différentiel de Wey-Hörmander associé à la métrique $ \frac{dx^2}{1+x^2} + \frac{d\xi^2}{1+\xi^2} $.
\\\\La classe des symboles $ S(\mu,m) $ associée à la métrique précédente est l'espace des fonctions régulières sur $ \R^d * \R ^d  $ qui vérifient $ | \partial^\alpha_x  \partial_\xi^\beta a (x , \xi )  | \leq C_{ \alpha, \beta  } <x>^{\mu-\alpha} < \xi >^{m-\beta }  $.
\\\\Ainsi, nous avons (voir \cite{Hormander} section 18.5, \cite{Robert} ou \cite{Bouclet}) que si $ a_1 \in S(\mu_1,m_1 ) $ et $ a_2 \in S( \mu_2, m_2 ) $ alors le commutateur $ [ Op(a_1),Op(a_2) ] $ est un opérateur pseudo-différentiel avec un symbole dans la classe $ S ( \mu_1+\mu_2-1  , m_2+m_2-1) $.
\\\\Par conséquent,
\begin{equation*}
[\sqrt H^{s_1+s_2},<x>^{-s_1}] \in S (s_2-1 , s_1+s_2-1) \subset S ( 0 , s_1 + s_2-1 ).
\end{equation*}
De plus, comme rappelé dans \cite{Martinez}, si $ q \in S(0,\mu ) $ alors pour tout $ s \in \R $, il existe une constante $ C> 0 $ telle que
\begin{equation*}
||Op(q) u ||_{H^{s-\mu}(\R^d)} \leq C || u ||_{H^s(\R^d)}.
\end{equation*}
Ainsi, nous pouvons prendre $ s= \mu =s_1+ s_2 -1 $ pour obtenir que
\begin{equation*}
|| [\sqrt H^{s_1+s_2},<x>^{-s_1}] u || _{L^2(\R^d) } \leq C ||u||_{H ^{s_1+s_2 -1}(\R^d) } \leq ||u||_{L^2 (\R^d)}.
\end{equation*}
De manière similaire, 
\begin{align*}
[ \  | \nabla | ^{s_1};<x>^{-s_2}] \in S (-s_2-1,s_1-1) \subset S (0,s_1-1),
\end{align*}
et nous pouvons conclure de la même façon pour le second point.
\\\\Pour le dernier point, nous avons
\begin{equation*}
[  \sqrt{H}^{s_1};<x>^{-s_2} ]  \sqrt{H}^{s_2-s_1} \in S (-1,s_2-1) \subset S (0,s_2-1).
\end{equation*}
Puis 
\begin{equation*}
|| [  \sqrt{H}^{s_1};<x>^{-s_2} ]  \sqrt{H}^{s_2-s_1} u ||_{L^2(\R^d)} \leq C ||u||_{L^2(\R^d)},
\end{equation*}
et il suffit de remplacer $ u $ par $ \sqrt{H}^{s_1-s_2} u $ pour obtenir le résultat désiré. \hfill $ \boxtimes $
\begin{lem} \label{2gradient} Soit $ s \geq 0  $ alors il existe deux constantes $ C_1,C_2 > 0 $ telles que pour tout $ f \in C^\infty_0(\R^d) $,
\begin{equation*}
 C_1 \times || \ |\nabla| ^s(  f   ) ||_{L^2(\R^d)} \leq || \nabla^s(  f   ) ||_{L^2(\R^d)} \leq C_2 \times || \ |\nabla| ^s(  f   ) ||_{L^2(\R^d)}.
\end{equation*}
\end{lem}
\textit{Preuve.} En utilisant l'égalité de Plancherel, il suffit de remarquer que la fonction 
\begin{equation*}
b ( \xi_1, \ . \ . \  . \ , \xi_d ) = \frac{   \left(  \underset{i}{\sum} \xi_i^2 \right) ^{s/2}  }{ \sqrt{ \underset{i}{\sum} \xi_i^{2s}   }  },
\end{equation*}
est positive, continue sur $ \R^d \setminus \lbrace 0  \rbrace $ et homogène (c'est à dire que pour tout $ \xi \in \R^d$ et $ \lambda \in \R^*$, $ b( \lambda \xi ) =  b ( \xi) $). \hfill $ \boxtimes $
\\\\Ces différents lemmes établis, nous pouvons passer à la preuve de l'effet régularisant.
\subsection{Preuve de (\ref{2effectregularisant1})} 
\textbf{\'Etape 1:}
\\\`A l'aide du calcul pseudo différentiel, soit le théorème 2.6.5 de \cite{Martinez}, on trouve
\begin{align*}
& \left[ \frac{x.D_x}{<x>^\alpha} ; \Delta  \right] 
 \\ = \ &  \frac{x.D_x}{<x>^\alpha} \Delta-\Delta  \frac{x.D_x}{<x>^\alpha} 
\\ = \ &  Op \left( \frac{x. \xi}{ <x>^\alpha} \right) Op \left( -\xi^2 \right) - Op \left(  -\xi^2 \right) Op \left( \frac{x.\xi }{<x>^\alpha} \right) 
\\  = \ & Op \left(  -2i  \times \left( \frac{\xi^2}{<x>^\alpha } -\alpha \frac{(x. \xi)^2}{ <x>^{\alpha+2} } \right) + 2 \alpha (d+2) \frac{ x. \xi}{<x>^{\alpha+2}} +  2\alpha ( \alpha+2 ) \frac{ x. \xi \  x^2 }{<x>^{\alpha+4}} \right).
\end{align*}
Puis, en utilisant que $ \alpha < 1 $, on obtient
\begin{align*}
& \Re  \left( i \int_{\R^d } Op \left( -2i \times  \left( \frac{\xi^2}{<x>^\alpha } -\alpha \frac{(x. \xi)^2}{ <x>^{\alpha+2} } \right)  u (x) \ \ \times \overline{u} (x) \ dx \right) \right)
 \\\\  = \ & 2 \Re \left(  \int_{\R^d}   \frac{- \Delta u }{<x>^\alpha}  \ \  \overline{u} - \alpha \frac{(x.D_x)^2u + i  ( x.D_x ) u}{<x>^{\alpha+2}} \ \  \overline{u}   \  dx  \right)
 \\\\  = \ & 2 \Re \left(  \int_{\R^d}   \nabla u .  \nabla ( \frac{ \overline{u} }{<x>^{\alpha+2}} ) -\alpha ( x.\nabla u )  \times \   div \left( \frac{ x  \overline{u} }{<x>^{\alpha+2}}  \right) - i \alpha  \frac{(x. D_x ) u }{<x>^{\alpha+2}} \ \ \overline{u} \ dx  \right)
 \\\\  = \ & 2 \int_{\R^d } \left( \frac{|\nabla u |^2}{<x>^\alpha}   -\alpha \frac{( x.\nabla u )^2}{<x>^{\alpha+2}} \ dx \right) 
  \\ & \hspace*{3cm} +2 \alpha \ \Re \left(  \int_{\R^d}  (\alpha+2) \frac{x^2 ( x.\nabla u ) }{<x>^{\alpha+4}}  \ \  \overline{u} - (d+1) \frac{ x.\nabla  u    }{<x >^{\alpha+2}} \ \ \overline{u}  \ dx \right)
 \\\\ \geq \ & 2(1-\alpha) \times \bigg| \bigg| \frac{1}{<x>^{\alpha/2}} \nabla u \bigg| \bigg|^2_{L^2(\R^d)} 
 \\ & \hspace*{3cm} +2 \alpha \ \Re \left(  \int_{\R^d}  (\alpha+2) \frac{x^2 ( x.\nabla u ) }{<x>^{\alpha+4}}  \ \  \overline{u} - (d+1) \frac{ x.\nabla  u    }{<x >^{\alpha+2}} \ \ \overline{u}  \ dx \right).
 \end{align*}
Grâce au lemme \ref{2inter}, on établit
\begin{align*} \Re \bigg( i & \left. \int_{\R^d} Op \left( 2 \alpha (d+2) \frac{ x. \xi}{<x>^{\alpha+2}} +  2\alpha ( \alpha+2 ) \frac{ x. \xi \  x^2 }{<x>^{\alpha+4}} \right) u \ \  \overline{u}  \ dx \right)
\\ & + 2 \alpha \ \Re \left(  \int_{\R^d}  (\alpha+2) \frac{x^2 ( x.\nabla u ) }{<x>^{\alpha+4}}  \ \  \overline{u} - (d+1) \frac{ x.\nabla  u    }{<x >^{\alpha+2}} \ \ \overline{u}  \ dx \right) \leq C ||u||^2_{H^{1/2}(\R^d)}.
\end{align*}
Ainsi, pour tout $ \alpha < 1 $, il existe une constante $ C> 0 $ telle que pour tout $ u \in H^{1/2}(\R^d) $,
\begin{equation*}
\Re   \bigg( i \int_{\R^d }  \left[ \frac{x.D_x}{<x>^\alpha} ,  \Delta  \right] u (x) \overline{u} (x) dx \bigg) \geq 2(1-\alpha) \bigg| \bigg| \frac{1}{<x>^{\alpha/2}} \nabla u \bigg| \bigg|_{L^2(\R^d)}^2 -C ||u||_{H^{1/2}(\R^d)}^2.
\end{equation*}
\\De manière similaire, on a  
\begin{equation*}
  \left[  \frac{x.D_x}{<x>^\alpha} ; -x^2 \right] = Op \left( \frac{2ix^2}{<x>^\alpha} \right),
 \end{equation*}
et donc
\begin{align*}
\Re \bigg( i \int_{\R^d}  \left[ \frac{x.D_x}{<x>^\alpha} ; -x^2 \right] u(x) \ \ \overline{u(x)} \bigg) \ dx & = -2 \int_{\R^d} \frac{x^2}{<x>^\alpha} |u(x)|^2 \  dx 
\\ & \geq -C ||u||^2_{   \overline{H}^{(2-\alpha)/2}(\R^d)}.
\end{align*}
Finalement, nous avons montré que pour tout $ \alpha < 1 $, il existe une constante $ C> 0 $ telle que pour tout $ u \in \overline{H}^{  (2-\alpha)/2 }(\R^d) $
\begin{equation} \label{2fin1}
\Re \bigg( i \int_{\R^d }  \left[ \frac{x.D_x}{<x>^\alpha} ,  -H   \right] u (x) \overline{u} (x) dx \bigg) \geq 2(1-\alpha) \bigg| \bigg| \frac{1}{<x>^{\alpha/2}} \nabla u \bigg| \bigg|_{L^2(\R^d)}^2 -C ||u||_{  \overline{H}^{(2-\alpha)/2}(\R^d)}^2.
\end{equation}
\textbf{\'Etape 2:} 
\\Si nous choisissons  $ u = e^{-itH} u_0 $ alors nous trouvons
\begin{align*}
& -i \int_{\R^d }  \left[ \frac{x.D_x}{<x>^\alpha} ,  H  \right] u (t,x) \ \ \overline{u} (t,x)  \ dx 
\\  = \ & - i \int_{\R^d} \frac{x.\nabla \partial_t u(t,x)} {<x>^\alpha} \ \ \overline{ u(t,x)}  \ dx +  i\int_{\R^d} \frac{x.D_x}{<x>^\alpha} u(t,x)  \ \ \overline{Hu(t,x)} \ dx
\\   = \ & - i \int_{\R^d} \frac{x.\nabla_x \partial_t u(t,x) }{<x>^\alpha} \ \ \overline{ u(t,x)}  \ dx -i  \int_{\R^d} \frac{x.\nabla_x}{<x>^\alpha} u(t,x) \ \ \overline{\partial_t u(t,x)} \ dx
\\ = \ & -i \partial_t \left( \int_{\R^d }  \frac{x. \nabla_x}{<x>^\alpha} u (t,x) \ \ \overline{u(t,x)}  \ dx  \right) .
\end{align*}
Ainsi, grâce à (\ref{2fin1}), on obtient pour $ T \geq 0 $,
\begin{align*}
&  2(1-\alpha) \int_0^T \bigg| \bigg| \frac{1}{<x>^{\alpha/2}} \nabla( e^{-itH } u_0 ) \bigg| \bigg|^2_{L^2(\R^d)} \ dt &
 \\  \leq \ & CT||u_0||^2_{H^{(2-\alpha)/2}(\R^d)}  + \Re \left( i \int_{\R^d}   \frac{x.\nabla_x  u_0 }{<x>^\alpha}  \ \ \overline{ u_0} - \frac{x.\nabla_x e^{-iTH} u_0 }{<x>^\alpha} \ \ \overline{e^{-iTH} u_0}   \ dx  \right) . 
 \end{align*} 
Puis, par le lemme \ref{2interbis}, on trouve pour tout $ \alpha \in  ] 0 ,1 [ $, l'existence d'une constante $ C > 0 $ telle que pour tout $ T \geq 0 $ et $ u_0 \in \overline{H}^{(2-\alpha)/2}(\R^d)$,
\begin{equation*}
\int_0^T \bigg| \bigg| \frac{1}{<x>^{\alpha/2 }} \nabla e^{-itH} u_0 \bigg| \bigg|^2_{L^2(\R^d)} \ dt \leq   CT ||u_0||_{  \overline{H}^{(2-\alpha)/2}(\R^d)}^2.
\end{equation*}
\textbf{\'Etape 3 :}
\\On prend $ \alpha = 1 - 2 \epsilon $ avec $ \epsilon \in ] 0 , \frac{1}{2} [ $ pour avoir
\begin{equation*}
\int_0^T \bigg| \bigg| \frac{1}{<x>^{1/2-\epsilon}} \nabla e^{-itH} u_0 \bigg| \bigg|^2_{L^2(\R^d)} \ dt \leq   CT ||u_0||_{ \overline{H}^{1/2+\epsilon}(\R^d)}^2.
\end{equation*}
En utilisant le lemme \ref{2commutateur}, on obtient
\begin{align*} 
 & \int_0^T \bigg| \bigg| \frac{1}{<x>^{1/2-\epsilon}} H^{1/2-\epsilon/2} e^{itH}u_0 \bigg| \bigg|_{L^2(\R^d)}^2 \ dt   
 \\  \leq \ & \int_0^T \bigg| \bigg|  H^{1/2-\epsilon/2} \frac{1}{<x>^{1/2-\epsilon}} e^{itH}u_0 \bigg| \bigg|_{L^2(\R^d)}^2
 \\ & \hspace*{5cm} +\int_0^T \bigg| \bigg| \ \left[ \frac{1}{<x>^{1/2-\epsilon}} ; H^{1/2-\epsilon/2} \right] e^{itH}u_0 \bigg| \bigg|_{L^2(\R^d)}^2
 \\ \leq \ & \int_0^T \bigg| \bigg|  H^{1/2-\epsilon/2} \frac{1}{<x>^{1/2-\epsilon}} e^{itH}u_0 \bigg| \bigg| _{L^2(\R^d)}^2 + T ||u_0||^2_ {L^2(\R^d)}.
\end{align*} 
Puis, en utilisant la proposition \ref{comparaison}, on trouve
\begin{align*}
& \int_0^T \bigg| \bigg|  H^{1/2-\epsilon/2} \frac{1}{<x>^{1/2-\epsilon}} e^{itH}u_0 \bigg| \bigg|_{L^2(\R^d)}^2
\\ \leq \ & \int_0^T \bigg| \bigg| <x>^{1/2 + \epsilon/2} e^{itH} u_0 \bigg| \bigg| ^2_{L^2(\R^d)} dt + \int_0^T \bigg| \bigg| \nabla^{1 - \epsilon } \left( \frac{1}{<x>^{1/2-\epsilon}}e^{itH} u_0 \right) \bigg| \bigg|^2_{L^2(\R^d)} \ dt 
 \\  \leq \ & CT ||u_0||^2_{\overline{H}^{1/2+\epsilon}(\R^d)} + \int_0^T \bigg| \bigg| \nabla \left( \frac{1}{<x>^{1/2-\epsilon}}e^{itH} u_0 \right) \bigg| \bigg|^2_{L^2(\R^d)} \ dt 
 \\  \leq \ & CT ||u_0||^2_{ \overline{H}^{1/2+\epsilon}(\R^d)} + \int_0^T \bigg| \bigg|  \frac{1}{<x>^{1/2-\epsilon}} \nabla  \left( e^{itH} u_0 \right) \bigg| \bigg|^2_{L^2(\R^d)} \ dt 
 \\ \leq \ & CT ||u_0||^2_{ \overline{H}^{1/2+\epsilon}(\R^d)}.
\end{align*}
Par conséquent, nous trouvons 
\begin{equation*}
\int_0^T \bigg| \bigg| \frac{1}{<x>^{1/2-\epsilon}} H^{1/2-\epsilon/2} e^{-itH} u_0 \bigg| \bigg|^2_{L^2(\R^d)} \ dt   \leq  CT ||u_0||^2_{ \overline{H}^{1/2+\epsilon}(\R^d)}.
\end{equation*}
Et nous pouvons remplacer $ u_0 $ par $ H^{-1/4-\epsilon/2} u_0 $ pour prouver le théorème. \hfill $ \boxtimes $
\subsection{Preuve de (\ref{2effectregularisant2})} 
En utilisant le lemme \ref{2commutateur} et (\ref{2effectregularisant1}), on obtient
\begin{align*}
& \bigg| \bigg| \sqrt H ^{d/2-2\epsilon} \frac{1}{<x>^{1/2-\epsilon} } e^{itH} u_0 \bigg| \bigg|_{L^2([-2 \pi,2\pi]*\R^d)} 
\\  \leq \ & \bigg| \bigg|  \frac{1}{<x>^{1/2-\epsilon} } \sqrt H ^{d/2-2\epsilon} e^{itH} u_0 \bigg| \bigg|_{L^2([-2 \pi,2\pi]*\R^d)} 
\\ & \hspace*{3cm} + \bigg| \bigg| \ \left[ \sqrt H ^{d/2-2\epsilon} ; \frac{1}{<x>^{1/2-\epsilon} } \right]  e^{itH} u_0 \bigg| \bigg| _{L^2([-2 \pi,2\pi]*\R^d)}
\\  \leq \ & C ||u_0||_{   \overline{H}^{  \frac{d-1}{2} } (\R^d)}.
\end{align*}
Puis, en utilisant la proposition \ref{comparaison}, on établit
\begin{equation*}
\bigg| \bigg| \nabla ^{d/2-2\epsilon}  \left(  \frac{1}{<x>^{1/2-\epsilon} } e^{itH} u_0  \right) \bigg| \bigg|_{L^2([-2 \pi,2\pi]*\R^d)}  \leq C||u_0||_{   \overline{H}^{  \frac{d-1}{2} } (\R^d)}.
\end{equation*}
Et finalement, en utilisant les lemmes \ref{2commutateur} et \ref{2gradient}, nous pouvons conclure que 
\begin{align*}
& \bigg| \bigg|  \frac{1}{<x>^{1/2-\epsilon} } |\nabla | ^{d/2-2\epsilon} ( e^{itH} u_0 ) \bigg| \bigg|_{L^2([-2 \pi,2\pi]*\R^d)} 
\\   \leq \ & \bigg| \bigg| \ \left[ \frac{1}{<x>^{1/2-\epsilon} } ; | \nabla | ^{d/2-2\epsilon}  \right]  e^{itH} u_0 \bigg| \bigg|_{L^2([-2 \pi,2\pi]*\R^d)}  
\\ & \hspace*{3cm} + \bigg| \bigg|  \nabla ^{d/2-2\epsilon} \left( \frac{1}{<x>^{1/2-\epsilon} } e^{itH} u_0 \right)  \bigg| \bigg|_{L^2([-2 \pi,2\pi]*\R^d)} 
\\   \leq \ & C ||u_0||_{  \overline{H}^{  \frac{d-1}{2} } (\R^d)}. \hspace*{8cm} \boxtimes
\end{align*}
\section{Données initiales aléatoires et espaces de Sobolev}
De manière analogue à la section 4 de \cite{moi}, on démontre que la donnée initiale rendue aléatoire ne permet pas de gagner de dérivée dans $ L^2( \mathds{R}^d )  $.
\begin{thm} \label{2sobolev2} Sous les hypothèses (\ref{2hypothese1}), (\ref{2hypothese2}) et (\ref{2hypothese4}), pour tout $ s \geq 0 $,
\begin{equation*}
\mbox{si } u_0  \notin \overline{H}^s(\R^d) \mbox{ alors }  u^\omega_0 \notin H^s(\R^d)   \  \omega \ ps.
\end{equation*}
\end{thm}
Pour établir ce résultat, en analogie au théorème 52 de \cite{moi}, nous devons montrer le même type d'estimation que la proposition 30 de  \cite{moi} pour des fonctions propres quelconques de l'oscillateur harmonique. Cela justifie la proposition suivante :
\begin{prop} Pour tout $ s \geq 0 $, il existe deux constantes $ C_1,C_2 > 0 $ telles que pour tout $ n \in \mathds{N}  $,
\begin{equation} \label{2propre2}
C_1 \lambda_n ^ s \leq ||    \nabla ^s  h_n ||_{L^2(\R^d)}  \leq C_2\lambda_n ^ s.
\end{equation}
\end{prop}
\textit{Preuve.} Nous posons $ h = \frac{1}{  \lambda_n^2 } $ et $ \Phi_h( x) =  \frac{1}{h^{d/4}}  \times h_n( \lambda_n x ) $ pour que $ ( -h^2 \Delta + x^2 -1 )  \Phi_h = 0 $ et $  || \Phi_h ||_{L^2(\R^d)}  =1 $.
\\ Pour démontrer (\ref{2propre2}), il suffit d'établir qu'il existe une constante $ C_1> 0 $ telle que pour tout $ h > 0 $, 
\begin{equation*}
h^s || \nabla^s  \Phi_h  ||_{L^2(\R^d)} \geq C_1.
\end{equation*}
Raisonnons par l'absurde et supposons que 
\begin{equation} \label{2absurdite}
\underset{h \rightarrow 0}{\underline{\lim}} \ h^s || \nabla^s  \Phi  ||_{L^2(\R^d)} = 0.
\end{equation}
D'après le théorème 2 de \cite{burq8}, il existe une mesure positive $ \mu \in  \mathcal{M}_+ ( \R^d * \R^d ) $ telle que pour toute fonction $ a \in C_0^\infty (\R^d * \R^d   )$ ,
\begin{equation*}
\lim_{ h \rightarrow 0 } < a ( x, hD_x )  \Phi_h , \Phi_h > _{ L^2(\R^d) * L ^2(\R^d)   } = \int_{ \R^d* \R^d  } tr ( a(x, \xi ) )  \ \mu(  d x d \xi  ).
\end{equation*}
Rappelons la définition suivante:
\begin{dfn}
On dit que $ (x, \xi) \in Supp(  \mu )^c $ si et seulement si il existe $ r > 0 $ tel que pour tout $ \phi \in C_0^\infty ( B(x,r) \times B(\xi,r)  )$,
\begin{equation*}
\int_{ \R^d * \R^d   } \phi (x,\xi) \ \mu(dx,d\xi) = 0.
\end{equation*}
\end{dfn}
De manière similaire à la proposition 40 de \cite{moi}, si $ a \in C_0^\infty (\R^d * \R^d  ) $  avec $ Supp ( a)  \cap  \lbrace (x,\xi) / x^2 + \xi^2  = 1   \rbrace =  \varnothing $ alors pour tout $ N \in \N $, il existe $ E_N \in Op ( T^{-2} ) $ et $ R_N \in Op ( T^{-(N+1)} ) $ tels que
\begin{equation*}
E_N \circ (-h^2 \Delta + |x|^2 - 1  )   = a ( x, hD_x ) - h^{N+1} R_N.
\end{equation*}
Par conséquent
\begin{equation*}
< a ( x, hD_x )  \Phi_h , \Phi_h > _{ L^2(\R^d) * L ^2(\R^d)   } =  h^{N+1} < R_N \Phi_h , \Phi_h > _{ L^2(\R^d) * L ^2(\R^d)   },
\end{equation*}
puis
\begin{equation*}
 \int_{ \R^d* \R^d  }  a(x, \xi )   \mu(  d x d \xi  ) = 0.
\end{equation*}
Et finalement, nous établissons que
\begin{equation*}
Supp ( \mu ) \subset \lbrace  (x,\xi) /  x^2  + \xi^2 = 1 \rbrace.
\end{equation*}
Toujours de manière similaire à la proposition 40 de \cite{moi}, si $ a \in C_0^\infty (\R^d * \R^d  ) $  avec $ Supp ( a)  \cap  \lbrace (x,\xi) / \xi^2  = 0   \rbrace =  \varnothing $ alors pour tout $ N \in \N $, il existe $ E_N \in Op ( S^{-s} ) $ et $ R_N \in Op ( S^{-(N+1)} ) $ tels que
\begin{equation*}
E_N \circ \sum_{i=1}^d | h D_{x_i} | ^{s}   = a ( x, hD_x ) - h^{N+1} R_N.
\end{equation*}
Or d'après \cite{Martinez} et (\ref{2absurdite}), on trouve
\begin{align*}
 \lim_{ h \rightarrow 0 }  | < E_N \circ \sum_{i=1}^d | h D_{x_i} | ^{s} \Phi_h , \Phi_h > _{ L^2(\R^d) * L ^2(\R^d)   } | & \leq \lim_{ h \rightarrow 0 } || E_N \circ \sum_{i=1}^d | h D_{x_i} | ^{s} \Phi_h ||_{L ^2(\R^d)}
 \\ & \leq \lim_{ h \rightarrow 0 } || \sum_{i=1}^d | h D_{x_i}| ^{s} \Phi_h ||_{L ^2(\R^d)} = 0
\end{align*}
Par conséquent
\begin{equation*}
 \int_{ \R^d* \R^d  }  a(x, \xi )   \mu(  d x d \xi  ) = 0,
\end{equation*}
et nous établissons que
\begin{equation*}
Supp ( \mu ) \subset \lbrace  (x,\xi) / \xi^2 = 0 \rbrace.
\end{equation*}
Ensuite, pour $ a \in  C_0^\infty(  \R^d * \R^d ) $ alors
\begin{align*}
0 & = \int_ {\R^d} [   -h^2 \Delta + |x|^2 - 1 ; h^{-1} Op_h(a)   ] \Phi_h \overline{\Phi_h}
\\ & = \frac{1}{i} \int_{\R^d} \lbrace -h^2 \Delta + |x|^2 - 1 ; Op_h(a) \rbrace \Phi_h \overline{\Phi_h} + h \times  \int_ {\R^d} Op_h(R ) \Phi_h \overline{\Phi_h}.
 \end{align*}
Ainsi, nous déduisons que pour toute fonction $ a \in C^\infty_0 ( \R^d) $,
\begin{equation} \label{2support mesure invariant}
\int_{\R^d*\R^d}  \left( \xi \partial_x a - x \partial_\xi a      \right) \  d\mu(x,\xi) = 0.
\end{equation}
Soit alors $ (x,\xi ) \in \R^d * \R^d $ et posons, pour $ t \in \R $, 
\begin{equation*}
\left\{
    \begin{array}{ll}
        & x(t) = x \cos(t) + \xi \sin(t) ,
        \\ & \xi(t) = \xi \cos(t) - x \sin(t). 
    \end{array}
\right.
\end{equation*}
c'est à dire 
\begin{equation*}
  \left\{
      \begin{aligned}
       &  \dot{x}(t)= \xi (t) \mbox{ avec } x(0)=x,\\
       &  \dot{\xi}(t)=-x(t) \mbox{ avec } \xi(0)= \xi. \\
      \end{aligned}
    \right.
\end{equation*}
D'après (\ref{2support mesure invariant}), on obtient pour tout $ t \in \R $ et $ a \in C_0^\infty(\R^d*\R^d) $,
\begin{equation*}
\int_{\R^d*\R^d}  a( x \cos(t) + \xi \sin(t) ,\xi \cos(t) - x \sin(t)) \ d\mu(x,\xi)  =\int_{\R^d*\R^d}  a(x,\xi) \ d\mu(x,\xi).
\end{equation*}
Par conséquent, si $ (x_0,\xi_0) \in Supp( \mu ) $ alors pour tout $ r > 0 $, il existe $ a \in C^\infty_0 ( B ( (x_0,\xi_0) ,r )) $ telle que pour tout $ t \in \R $,
\begin{equation*}
\int_{\R^d*\R^d}  a( x \cos(t) + \xi \sin(t) ,\xi \cos(t) - x \sin(t))  \   d\mu(x,\xi)  \neq 0. 
\end{equation*}
Mais
\begin{align*}
 Supp & \left( a( x \cos(t) + \xi \sin(t) ,\xi \cos(t) - x \sin(t))     \right) 
 \\ & \subset \left\{  (x,\xi) \in \R^d* \R^d /  ( x \cos(t) + \xi \sin(t)  , \xi \cos(t) - x \sin(t )  ) \in B ( (x_0,\xi_0) ,r ) \right\} 
 \\ & \subset  B(  \cos(t) x_0 - \sin(t) \xi_0  ,  2 r ) \times B( \sin(t) x_0 + \cos(t) \xi_0  ,  2 r  ),
 \end{align*}
et donc, pour tout $ t\in \R$,  $ ( \cos(t) x_0 - \sin(t) \xi_0 , \sin(t) x_0 + \cos(t) \xi_0  )  \in Supp( \mu ) $.
\\\\Mais pour $ \xi_0 = 0, x_0^2 = 1 $ alors $ \sin(t) x_0 + \cos(t) \xi_0 = \sin(t) x_0 = 0 $ est impossible et donc la proposition est démontrée par l'absurde. \hfill $ \boxtimes $
\\\\Ensuite, pour une fonction $ \chi \in C_0^\infty (\R^d)$ telle que  $\chi(x)=1 $ si $ |x| \leq 1 $, $ \chi(x) = 0  $ si $ |x| \geq 2 $ et $ 0 \leq \chi \leq 1 $, définissons
\begin{align*}
& \sigma_N^2 = \sum_{n \in \mathds{N}} \chi^2  \left( \frac{ \lambda_n^2 }{N^2} \right) |c_n|^2 \lambda_n^{2s}  \underset{ N \rightarrow \infty }{\longrightarrow} \infty,
\\ & S_N = || \chi \left(  \frac{H}{N^2}  \right) u_0  ||_{H^s(\R^d)},
\\ & M = \displaystyle{  \sup_{N \in \N ^*} } S_N.
\end{align*}
Passons à la preuve du théorème \ref{2sobolev2}. En analogie à la preuve du théorème 52 de \cite{moi}, il suffit d'établir que 
\begin{equation*}
P( M = \infty ) > 0.
\end{equation*}
Grâce à (\ref{2propre2}) et aux hypothèses (\ref{2hypothese2bis}) et (\ref{2hypothese4}), on trouve
\begin{align*}
& E \left(  || \chi \left( \frac{H}{N^2}  \right)  u_0 ||^2_{H^s(\R^d)} \right) 
\\ \geq \ & E \left( \sum_{n,m} \chi \left( \frac{\lambda_n^2}{N^2} \right) \chi \left( \frac{\lambda_m^2}{N^2} \right) c_n \overline{c_m} \ g_n( \omega ) \overline{g_m ( \omega )} \ \int_{\R^d} \nabla^s ( h_n )  \nabla^s( h_m ) \ dx   \right)
\\  \geq  \ & E \left( \sum_{n \in \mathds{N} } \chi^2 \left( \frac{\lambda_n^2}{N^2} \right) | c_n |^2 | g_n \omega )|^2  || \nabla^s ( h_n )  ||^2_{L^2(\R^d)}   \right)
\\ \geq \ & C_1 \sigma_N^2.
\end{align*}
Par conséquent, grâce à l'inégalité de Zygmound, soit le lemme 53 de \cite{moi} , on établit
\begin{align*}
P \left( M^2 \geq \frac{C_1 \sigma_N^2}{2} \right)  \geq P \left( S_N^2 \geq  \frac{ C_1  \sigma_N^2}{2} \right) &  \geq P \left( S_N^2 \geq  \frac{ || \chi \left( \frac{H}{N^2}  \right)  u_0 ||^2_{H^s(\R^d)}  }{2} \right)
 \\ &  \geq \frac{1}{4} \times \frac{E \left( || \chi \left( \frac{H}{N^2}  \right)  u_0 ||^2_{H^s(\R^d)} \right) ^2}{E  \left( || \chi \left( \frac{H}{N^2}  \right)  u_0 ||^4_{H^s(\R^d)}\right) }.
\end{align*}
Puis, grâce à (\ref{2propre2}), au lemme \ref{2moments} et l'hypothèse (\ref{2hypothese2bis}), on a 
\begin{align*}
& E \left( || \chi \left( \frac{H}{N^2}  \right)  u_0 ||^4_{H^s(\R^d)}  \right) 
\\  \leq \ &  E \left( \sum_{n, m}   \chi \left( \frac{\lambda_n^2}{N^2} \right) \chi \left( \frac{\lambda_m^2}{N^2} \right) c_n \overline{c_m} \ g_n( \omega ) \overline{g_m ( \omega )} \ \int_{\R^d} \nabla^s ( h_n )  \nabla^s( h_m ) \ dx  \right)^2
\\   & \hspace*{4cm} + E \left( \sum_{n,m}   \chi \left( \frac{\lambda_n^2}{N^2} \right) \chi \left( \frac{\lambda_m^2}{N^2} \right) c_n \overline{c_m} \ g_n( \omega ) \overline{g_m ( \omega )} \right)^2
\\ \leq \ & C E \left(  \sum_{n_1,n_2,n_3,n_4 } \chi \left( \frac{\lambda_{n_1}^2}{N^2} \right) \chi \left( \frac{\lambda_{n_2}^2}{N^2} \right) \chi \left( \frac{\lambda_{n_3}^2}{N^2} \right) \chi \left( \frac{\lambda_{n_4}^2}{N^2} \right) \times c_{n_1} \overline{c_{n_2}} c_{n_3} \overline{c_{n_4}} \right.
\\  & \hspace*{0.5cm}  \times || \nabla^s( h_{n_1} ) ||_{L^2( \R^3 )} || \nabla^s( h_{n_2} ) ||_{L^2( \R^3 )}|| \nabla^s( h_{n_3} ) ||_{L^2( \R^3 )}|| \nabla^s( h_{n_4} ) ||_{L^2( \R^3 )} \bigg)
\\   &  \hspace*{7cm} + C E \left( \sum_{ n}   \chi^2 \left( \frac{\lambda_n^2}{N^2} \right) | c_n |^2  \right)^2
\\ \leq \ & C_2 \sigma_N ^4.
\end{align*}
Par conséquent,
\begin{equation*}
P \left( M^2 \geq \frac{C_1 \sigma_N^2}{2}  \right) \geq \frac{1}{4} \times \frac{C_1^2}{C_2},
\end{equation*}
puis, en utilisant un théorème de convergence monotone, on trouve
\begin{equation*} 
P \left( M = \infty  \right) \geq \frac{1}{4} \times \frac{C_1^2}{C_2}.
\end{equation*}
Et le théorème est démontré. \hfill $ \boxtimes $
\section{L'argument de point fixe}
Dans cet partie, on établit des estimées qui seront utiles pour appliquer un théorème de point fixe de Picard. On commence par montrer deux lemmes préliminaires.
\begin{lem} \label{2injection}
Soient $ (q,r) \in [2,\infty [ \times [2,\infty ] , \ s,s_0 \geq 0 $ et supposons que $ s-s_0 >  \frac{d}{2} - \frac{2}{q} - \frac{d}{r} $, alors il existe deux constantes $ \kappa,C > 0 $ telles que pour tout $ T \geq 0 $ et $  u \in \overline{X}_T^s $,
 \begin{equation*}
|| u ||_{ L^q( [-T,T] , \overline{W}^{s_0,r}(\R^d)) }  \leq  C T ^\kappa ||u||_{  \overline{X}_T^s }.
\end{equation*}
\end{lem}
\textit{Preuve.} Soit $ \epsilon > 0  $ alors il existe $ \kappa_\epsilon > 0  $ tel que 
\begin{align*}
|| u ||_{ L^q( [-T,T] , \overline{W}^{s_0,r}(\R^d)) }  \leq T^{\kappa_\epsilon} || u ||_{ L^{q+\epsilon}( [-T,T] , \overline{W}^{s_0,r}(\R^d)) }.
\end{align*}
Or le couple $ (  q+\epsilon  ,  \frac{2d(q+\epsilon)}{dq+d\epsilon-4}    ) $ est admissible avec 
\begin{equation*}
   \overline{W}^{s,\frac{2d(q+\epsilon)}{dq+d\epsilon-4}}(\R^d)  \hookrightarrow \overline{W}^{s_0,r}(\R^d) \mbox{ si } s-s_0 \geq  \frac{d}{2} - \frac{2}{q+\epsilon} - \frac{d}{r}.
\end{equation*}
Mais, comme $ s-s_0 >  \frac{d}{2} - \frac{2}{q} - \frac{d}{r} $ alors il existe $ 0 < \epsilon \ll 1 $ tel que $ s-s_0 \geq  \frac{d}{2} - \frac{2}{q+\epsilon} - \frac{d}{r}$. \hfill $ \boxtimes $
\begin{lem} \label{2gradient2} Soit $ s \geq 0 $ alors il existe une constante $ C > 0 $ telle que pour toutes fonctions f et g dans $ \mathcal{S}   (\R^d) $,
\begin{equation*}
|| \nabla^s(  f g  ) ||_{L^2(\R^d)} \leq C \times \left( \ || \ | \nabla|^s(  f ) \times  g   ||_{L^2(\R^d)} + ||  f \times   | \nabla   | ^s (  g  ) ||_{L^2(\R^d)} \ \right) .
\end{equation*}
\end{lem}
\textit{Preuve.} Par la transformée de Fourier et le lemme \ref{2gradient}, on obtient
\begin{align*}
|| \nabla^s(  f g  ) ||_{L^2(\R^d)} & \leq C \times || \ |\xi|^s  \mathcal{F} (  f g   ) ||_{L^2(\R^d)}
\\  & \leq  C \times || \ |\xi|^s (  \mathcal{F} (  f ) *   \mathcal{F} ( g )  )  ||_{L^2(\R^d)}.
\end{align*}
Or pour tous $ \xi $ et $ \eta $ dans $ \R^d $ , 
\begin{align*}
|\xi|^s \leq ( |\eta| + |\xi-\eta| )^s \leq C_s \times ( |\eta|^s +  |\xi-\eta|^s  ).
\end{align*}
Ainsi,
\begin{align*}
|| \nabla^s(  f g  ) ||_{L^2(\R^d)} & \leq C_s \times \left( \ ||   ( |.|^s \mathcal{F} (  f ) ) *   \mathcal{F} ( g )      ||_{L^2(\R^d)} + ||  \mathcal{F} (  f ) *   ( |.|^s \mathcal{F} (  g ) )      ||_{L^2(\R^d)} \ \right) &
\\ & \leq C_s \times \left( \ ||   \mathcal{F} (  |\nabla|^s f \times  g )      ||_{L^2(\R^d)} + ||  \mathcal{F} (  f ) *   \mathcal{F} ( |\nabla|^s g )      ||_{L^2(\R^d)} \ \right) 
\\ & \leq C_s \times \left( \ ||  \   |\nabla|^s f \times  g       ||_{L^2(\R^d)} + ||    f  \times   |\nabla|^s g       ||_{L^2(\R^d)} \ \right) . & \boxtimes
\end{align*}
Puis, on établit les estimées attendues.
\begin{prop} \label{2estim1}
Soit $ s > \frac{d}{2} - \frac{2}{p-1} $ alors il existe deux constantes $ C >  0 $ et $ \kappa > 0 $ telles que si nous supposons
\begin{equation*}
 || e^{-itH} u_0 ||_{ L^p  (  [-2 \pi, 2 \pi ] ,  L^ \infty ( \R^d )  )  } \leq \lambda 
\end{equation*}
pour un certain $ \lambda $, alors pour tout $ 0 < T \leq 1 $, $ v \in \overline{X}^s_T $ et $ f_i = v $ ou $ f_i = e^{-itH} u_0 $,
\begin{equation*}
|| \ |\nabla| ^s( v  ) \times   \prod_{ i = 2 } ^p  f_i  ||_{  L^1( [-T,T]  , L ^2(\R^d)) } \leq C  T ^\kappa(  \lambda ^p + ||v||^p_{\overline{X}^s_T}  ),
\end{equation*}
et 
\begin{equation*}
|| <x>^s \times   v  \times   \prod_{ i = 2 } ^p  f_i  ||_{  L^1( [-T,T]  , L ^2(\R^d)) } \leq C  T ^\kappa(  \lambda ^p + ||v||^p_{\overline{X}^s_T}  ).
\end{equation*}
\end{prop}
\textit{Preuve.} D'après l'inégalité de Hölder et la proposition \ref{comparaison},
\begin{align*}
& || \ |\nabla|^s( v  ) \times   \prod_{ i = 2 } ^p  f_i  ||_{  L^1( [-T,T]  , L ^2(\R^d)) } 
\\  \leq \ & || \ | \nabla|^s( v  ) ||_{  L^\infty ( [-T,T], L^2(\R^d)) } \times \prod_{ i = 2 } ^p || f_i ||_{  L^{p-1} ( [-T,T], L^\infty(\R^d)) }
\\ \leq \ & C || v ||_{  L^\infty ( [-T,T], \overline{H}^s (\R^d)) } \times \prod_{ i = 2 } ^p || f_i ||_{  L^{p-1} ( [-T,T], L^\infty(\R^d)) },
\end{align*}
et
\begin{align*}
& || <x>^s \times  v   \times   \prod_{ i = 2 } ^p  f_i  ||_{  L^1( [-T,T]  , L ^2(\R^d)) } 
\\ \leq \ &  || <x>^s v   ||_{  L^\infty ( [-T,T], L^2(\R^d)) } \times \prod_{ i = 2 } ^p || f_i ||_{  L^{p-1} ( [-T,T], L^\infty(\R^d)) }
\\ \leq \ & || v ||_{  L^\infty ( [-T,T], \overline{H}^s (\R^d)) } \times \prod_{ i = 2 } ^p || f_i ||_{  L^{p-1} ( [-T,T], L^\infty(\R^d)) }.
\end{align*}
Si $ f_i = v $ alors comme $ s > \frac{d}{2} - \frac{2}{p-1} $, nous pouvons utiliser le lemme \ref{2injection} pour obtenir
\begin{equation*}
|| v ||_{  L^{p-1} ( [-T,T], L^\infty(\R^d)) }  \leq C T^{\kappa  } || v ||_{  \overline{X}^s_T }.
\end{equation*}
Si $ f_i = e^{-itH} u_0 $ alors d'après l'inégalité de Hölder, 
\begin{align*}
|| e^{-itH} u_0  ||_{  L^{p-1} ( [-T,T], L^\infty(\R^d)) } & \leq T^{   \frac{1}{p-1} - \frac{1}{p} } || e^{-itH} u_0 ||_{ L^p  (  [-2 \pi, 2 \pi ] ,  L^ \infty ( \R^d)  )  } &
\\ & \leq T^{   \frac{1}{p(p-1)} }  \lambda. & \hspace*{1cm} \boxtimes 
\end{align*}
\begin{prop} \label{2estim2}
Soit $  \frac{d}{2} > s > 0 $ alors il existe deux constantes $ C >  0 $ et $ \kappa > 0 $ telles que si nous supposons que
\begin{equation*}
 || e^{-itH} u_0 ||_{ L^p  (  [-2 \pi, 2 \pi ] ,   \overline{W}^{  \frac{1}{8} , \infty } ( \R^d )  )  } \leq \lambda, 
\end{equation*}
et
\begin{equation*}
 || u_0 ||_{   \overline{H}^{  \frac{d-1}{2} } ( \R^d   )  } \leq \lambda 
\end{equation*}
pour un certain $ \lambda $, alors pour tout $ 0 < T \leq 1 $,
\begin{equation*}
|| <x>^s \times \ ( e^{-itH } u_0 )^p ||_{  L^1( [-T,T]  , L ^2(\R^d)) } \leq C  T ^\kappa  \lambda ^p.
\end{equation*}
\end{prop}
\textit{Preuve.} D'après l'inégalité de Hölder et la proposition \ref{comparaison}, on obtient
\begin{align*}
 & || <x>^s \times \ ( e^{-itH } u_0 )^p ||_{  L^1( [-T,T]  , L ^2(\R^d)) } 
\\  \leq \ & || <x>^{ \frac{d}{2}  } \times \ ( e^{-itH } u_0 )^p ||_{  L^1( [-T,T]  , L ^2(\R^d)) }
\\  \leq \ &   || <x>^{ \frac{d-1}{2}  } \times \  e^{-itH } u_0  ||_{  L^\infty( [-T,T]  , L ^2(\R^d)) } \times  || <x>^{ \frac{1}{2(p-1)}  } \times \  e^{-itH } u_0  ||_{  L^{p-1}( [-T,T]  , L ^\infty(\R^d)) }  ^{p-1} &
\\  \leq \ & C T^{1/p}  ||  u_0  ||_{  \overline{H} ^{  \frac{d-1}{2}}(\R^d)) } \times  ||   e^{-itH } u_0  ||_{  L^{p}( [-T,T]  , \overline{W}^{ \frac{1}{8} , \infty}(\R^d))  }  ^{p-1}
\\   \leq \ & C  T^{1/p} \lambda ^p.  \hspace*{10cm}  \boxtimes 
\end{align*}
\begin{prop} \label{2estim3}
Il existe $ s \in  ]  \frac{d}{2}   - \frac{2}{p-1} ; \frac{d}{2} [ $, $ C >  0 $ et $ \kappa > 0 $ tels que si nous supposons que
\begin{align*}
 || e^{-itH} u_0 ||_{ L^{2p}  (  [-2 \pi, 2 \pi ] ,  \overline{W} ^ { \frac{1}{7}  ,  \infty }  ( \R^d )  )  } & \leq \lambda,
 \end{align*}
et
 \begin{align*}
 ||u_0||_{  \overline{H}^{ \frac{d-1}{2} } (\R^d)} & \leq \lambda  
\end{align*}
pour un certain $ \lambda $, alors pour tout $ 0 < T \leq 1 $, $ v \in \overline{X}^s_T $ et $ f_i = v $ ou $ f_i = e^{-itH} u_0 $,
\begin{equation*}
|| \ |\nabla|^s( e^{-itH} u_0   ) \times   \prod_{ i = 2 } ^p  f_i  ||_{  L^1( [-T,T]  , L ^2(\R^d)) } \leq C  T ^\kappa(  \lambda ^p + ||v||^p_{\overline{X}^s_T}  ).
\end{equation*}
\end{prop}
\textit{Preuve.} Pour tout $ \epsilon \in ] 0, \frac{1}{2} [ $, d'après l'inégalité de Hölder, on a 
\begin{align*}
& || \ |\nabla|^s( e^{-itH} u_0   ) \times   \prod_{ i = 2 } ^p  f_i  ||_{  L^1( [-T,T]  , L ^2(\R^d)) } 
\\  \leq \ & \bigg| \bigg|   \frac{1}{ <x>^{1/2-\epsilon} } |\nabla|^s( e^{-itH} u_0   ) \bigg| \bigg|_{  L^2 ( [-T,T], L^2 (\R^d)) } \prod_{ i = 2 } ^p ||  <x>^{\frac{1}{p-1} * ( \frac{1}{2}  - \epsilon )} f_i ||_{  L^{2(p-1)} ( [-T,T],  L^{ \infty }(\R^d)) }.
\end{align*}
Puis, nous choisissons $ s = \frac{d}{2} - 2\epsilon $ avec $ \epsilon \ll 1 $ pour obtenir en utilisant (\ref{2effectregularisant2}) que
\begin{align*}
& || \ |\nabla| ^s( e^{-itH} u_0   ) \times   \prod_{ i = 2 } ^p  f_i  ||_{  L^1( [-T,T]  , L ^2(\R^d)) }   \\ \leq \ & \lambda \times \prod_{ i = 2 } ^p || <x>^{\frac{1}{p-1} * ( \frac{1}{2}  - \epsilon )} f_i ||_{  L^{2(p-1)} ( [-T,T],  L^{  \infty }(\R^d)) }
 \\ \leq \ & \lambda \times \prod_{ i = 2 } ^p ||  f_i ||_{  L^{2(p-1)} ( [-T,T],   \overline{W} ^{ \frac{1}{p-1} * ( \frac{1}{2}  - \epsilon ) + \epsilon , \frac{d}{\epsilon} +1 }(\R^d)) },
\end{align*}
Si $ f_i = e^{-itH} u_0 $, par interpolation, nous pouvons trouver l'existence d'une constante $ \kappa > 0 $ telle que
\begin{align*}
& || e^{-itH} u_0  ||_{   L^{2(p-1)} ( [-T,T],   \overline{W} ^{ \frac{1}{p-1} * ( \frac{1}{2}  - \epsilon ) + \epsilon , \frac{d}{\epsilon} +1 }(\R^d))  } 
\\   \leq \ & C \times T ^\kappa \times  || e^{-itH} u_0  ||^\theta_{   L^{2p} ( [-T,T],   \overline{W} ^{ s_0 , \infty }(\R^d))  } \times || e^{-itH} u_0  ||^{1-\theta}_{   L^{\infty} ( [-T,T],   \overline{W} ^{ \frac{d-1}{2} , 2 }(\R^d))  }
\end{align*}
où $ \theta = \frac{d-\epsilon}{d+\epsilon} $ et $ s_0 = (  \frac{1-\theta}{\theta} )( \frac{d-1}{2}  )+ \frac{1}{\theta(p-1)} ( \frac{1}{2} - \epsilon  ) + \frac{\epsilon}{\theta} $.
\\\\Or $ || e^{-itH} u_0  ||_{   L^{\infty} ( [-T,T],   \overline{W} ^{ \frac{d-1}{2} , 2 }(\R^d))  } \leq || u_0  ||_{   \overline{H} ^{ \frac{d-1}{2} }(\R^d)  } \leq \lambda $, puis comme 
\begin{align*}
s_0 = \frac{1}{2(p-1)} + C \epsilon + o (\epsilon) \leq \frac{1}{7}
\end{align*}
alors $ || e^{-itH} u_0  ||_{   L^{2p} ( [-T,T],   \overline{W} ^{ s_0 , \infty }(\R^d))  } \leq \lambda $ \\et donc $ || e^{-itH} u_0  ||_{   L^{2(p-1)} ( [-T,T],   \overline{W} ^{ \frac{1}{p-1} * ( \frac{1}{2}  - \epsilon ) + \epsilon , \frac{d}{\epsilon} +1 }(\R^d))  } \leq \lambda $. 
\\\\Si $ f_i = v $, comme $ s - \frac{1}{p-1} * (  \frac{1}{2} - \epsilon ) > \frac{d}{2} - \frac{1}{p-1} - \frac{d\epsilon}{d+\epsilon} $ (si $ \epsilon \ll \frac{1}{2(p-2)} $) alors par le lemme \ref{2injection}, on trouve 
\\\\ \hspace*{2cm} $ || v  ||_{  L^{2(p-1)} ( [-T,T],   \overline{W} ^{ \frac{1}{p-1} * ( \frac{1}{2}  - \epsilon ) + \epsilon , \frac{d}{\epsilon} +1 }(\R^d))  }  \leq C  T^\kappa  ||v||_{\overline{X}^s_T } . \hfill \boxtimes $ 
\\\\Il est donc légitime d'introduire la définition suivante :  
\begin{dfn}
Soit $ \lambda   \geq 0 $ et définissons $ E_0( \lambda ) $ comme l'ensemble des fonctions $ u_0 \in \overline{H}^{  \frac{d-1}{2} } (\R^d) $ qui vérifient
\begin{align*}\left\{
    \begin{array}{ll}
      ||u_0||_{  \overline{H}^{  \frac{d-1}{2} } (\R^d) } & \leq \lambda,   
       \\  ||  e^{-itH} u_0 ||_{ L^{2p}( [-2\pi,2\pi] , \overline{W}^{  \frac{1}{7} , \infty } (\R^d))  } & \leq \lambda.
    \end{array} 
\right.
\end{align*}
\end{dfn} 
Puis, on peut établir les deux théorèmes principaux de cette partie.
\begin{thm} \label{2pointfixe1} Il existe $ s \in ] \frac{d}{2}   - \frac{2}{p-1}  ; \frac{d}{2} [  $, $ C > 0 $ et $ \kappa > 0 $ tels que si $ u_0 \in E_0(\lambda ) $ pour un certain $ \lambda > 0 $ alors pour tout $ v \in \overline{X}_T^s $ et $ 0 < T \leq 1 $,
 \begin{align*}
\bigg| \bigg|  \int_0^t   e^{-i(t-s)H}  K \cos ( 2s )^{ \frac{d(p-1)}{2} -2  } |  e^{-isH} u_0 + v |^{p-1} \times (   e^{-isH} u_0 + v ) \ ds   \bigg| \bigg|_{  \overline{X}_T^s }  
\\ \leq C \times T^\kappa \times ( \lambda ^p +  ||v||^p_{  \overline{X}_T^s } ).
\end{align*} 
\end{thm} 
\textit{Preuve.} En utilisant les propositions \ref{2Stricharz2} et \ref{comparaison}, on obtient
\begin{align*}
& \bigg| \bigg|  \int_0^t   e^{-i(t-s)H} K \cos ( 2s )^{  \frac{d(p-1)}{2} -2  }  |  e^{-isH} u_0 + v |^{p-1} \times (   e^{-isH} u_0 + v ) \ ds   \bigg| \bigg|_{  \overline{X}_T^s }  
\\  \leq   C & || \ K \cos ( 2s )^{  \frac{d(p-1)}{2} -2  } \times  |  e^{-isH} u_0 + v |^{p-1} \times (   e^{-isH} u_0 + v )   ||_{   L^1( [-T,T] , \overline{H}^s(\R^d) )   }  
\\ \leq  C & || \   |  e^{-isH} u_0 + v |^{p-1} \times (   e^{-isH} u_0 + v )   ||_{   L^1( [-T,T] , \overline{H}^s(\R^d) )   }  
\\ \leq  C  & || \  \nabla^s  \left( \  |  e^{-isH} u_0 + v |^{p-1} \times (   e^{-isH} u_0 + v ) \right)   ||_{   L^1( [-T,T] , L^2 (\R^d) )   }  
\\ & \hspace*{1cm} +  C || \ <x>^s \times   |  e^{-isH} u_0 + v |^{p-1} \times (   e^{-isH} u_0 + v )   ||_{   L^1( [-T,T] , L^2 (\R^d) )   }.   
\end{align*}
Puis, en utilisant le lemme \ref{2gradient2} et les propositions \ref{2estim1}, \ref{2estim2} et \ref{2estim3}, nous pouvons trouver une constante $ \kappa > 0  $ telle que pour tout $ u_0 \in E_0(\lambda) $, $ 0 < T \leq 1 $ et $ v  \in \overline{X}^s_T $,
\begin{align*}
|| \  \nabla^s  \left(  |  e^{-isH} u_0 + v |^{p-1} \times (   e^{-isH} u_0 + v )  \right)   ||_{   L^1( [-T,T] , L^2 (\R^d) )   }  \leq C T ^\kappa ( \lambda^p +   ||v||^p_{  \overline{X}_T^s }   ), 
\end{align*}
et 
\begin{align*}
|| \ <x>^s \times   |  e^{-isH} u_0 + v |^{p-1} \times (   e^{-isH} u_0 + v )   ||_{   L^1( [-T,T] , L^2 (\R^d) )   }    \leq C T ^\kappa ( \lambda^p +   ||v||^p_{  \overline{X}_T^s }   ). 
\end{align*}
\begin{flushright}
$ \boxtimes $ 
\end{flushright}
De manière similaire, on peut démontrer le théorème suivant :
\begin{thm} \label{2pointfixe2}
Il existe $ s \in ] \frac{d}{2}   - \frac{2}{p-1}  ; \frac{d}{2} [  $ (le même que dans le théorème précédent), $ C > 0 $ et $ \kappa > 0 $ tels que si $ u_0 \in E_0(\lambda ) $ pour un certain $ \lambda > 0 $ alors pour tout $ 0 < T \leq 1 $ et $ v_1,v_2 \in \overline{X}_T^s $,
 \begin{align*}
& \bigg| \bigg|  \int_0^t   e^{-i(t-s)H}  K \cos (2s)^{ \frac{d(p-1)}{2} -2  }  |  e^{-isH} u_0 + v_1 |^{p-1} \times (  e^{-isH} u_0 + v_1 ) \ ds  
\\ & \hspace*{2cm} -  \int_0^t   e^{-i(t-s)H} K \cos (2s)^{ \frac{d(p-1)}{2} -2  }  |  e^{-isH} u_0 + v_2 |^{p-1} \times (  e^{-isH} u_0 + v_2 ) \ ds  \bigg| \bigg|_{  \overline{X}_T^s }  
\\ & \leq C  T ^\kappa \times ||v_1-v_2||_{ \overline{X}^s_T } \times ( \lambda ^{p-1} +  ||v_1||^{p-1}_{  \overline{X}_T^s } + ||v_2||^{p-1}_{  \overline{X}_T^s } ).
\end{align*} 
\end{thm}
\section{Solutions globales pour l'équation (NLS)}
Dans cette partie, on applique un théorème de point fixe pour établir l'existence de solutions globales pour l'équation (\ref{2schrodinger}). Comme dans \cite{moi}, on introduit l'équation suivante :
\begin{equation} \label{2schrodingerH} \tag{NLSH}
  \left\{
      \begin{aligned}
         & i \frac{ \partial u }{ \partial t } - H u = K\cos (2 t)^{   \frac{d(p-1)}{2} -2 } \times |u|^{p-1} u,
       \\  &  u(0,x)  =u_0(x),
      \end{aligned}
    \right.
\end{equation}
où $ p \geq 5 $ désigne un entier impair et $ K \in \lbrace -1, 1 \rbrace $.
\begin{thm}  \label{2existencebis}
Il existe $ s \in ] \frac{d}{2}   - \frac{2}{p-1} ; \frac{d}{2} [  $, $ C > 0 $ et $ \delta > 0 $ tels que pour tout $ 0 < T \leq 1 $, si $ u_0 \in E_0(\lambda ) $ avec $ \lambda < C \times   T ^{ - \delta  } $ alors il existe une unique solution à l'équation (\ref{2schrodingerH}) sur $ [-T,T] $ dans l'espace $ e^{-itH}u_0 + B_{\overline{X}^s_T } (0,  \lambda   ) $.
\end{thm}
\textit{Preuve.} Définissons
\begin{equation*}
L ( v) = -i \int_0^t   e^{-i(t-s)H} K \cos (2s)^{  \frac{d(p-1)}{2} -2  }  |   e^{-isH} u_0 + v(s) |^{p-1} (  e^{-isH} u_0 + v(s) ) \ ds,
\end{equation*}
et remarquons que $ u = e^{-itH} u_0 + v $ est l'unique solution de (\ref{2schrodingerH}) sur $ [-T,T] $ dans l'espace $ e^{-itH}u_0 + B_{\overline{X}^s_T}(0 , R ) $ si et seulement si v est l'unique point fixe de L sur $ B_{\overline{X}^s_T}(0 , R ) $.
\\\\ Selon les propositions \ref{2pointfixe1} et \ref{2pointfixe2}, il existe deux constantes $ C > 0 $ et $ \kappa > 0 $ telles que 
\begin{align*}
& || L (v)  ||_{  \overline{ X } ^s_T  } \leq C T^\kappa ( \lambda ^p + || v||^p _{ \overline{ X } ^s_T }) 
\\ & || L (v_1)-L(v_2)  ||_{  \overline{ X } ^s_T  } \leq C T^\kappa || v_1-v_2|| _{ \overline{ X } ^s_T } ( \lambda ^{p-1} + || v_1||^{p-1}_{ \overline{ X } ^s_T }+ || v_2||^{p-1}_{ \overline{ X } ^s_T } ).
\end{align*}
Par conséquent, si $ \lambda <   (  \frac{1}{8 C  T^\kappa }  )^{  \frac{1}{p-1} }    $ alors L est une application contractante de $ B_{\overline{X}^{s}_T } (0, \lambda   ) $ et le théorème suit. \hfill $ \boxtimes $
\begin{thm}  \label{2existence1}
Il existe $  s \in ]  \frac{d}{2}   - \frac{2}{p-1} ; \frac{d}{2} [ $, $ C_1,C_2 > 0 $ tel que si $ u_0 \in E_0(\lambda ) $ avec $ \lambda < C_1 $ alors il existe une solution globale à (\ref{2schrodinger}) dans l'espace $ e^{it\Delta}u_0 + B_{X^s } (0, C_2 ) $.
\end{thm}
\textit{Preuve.} Soit u donnée par le théorème \ref{2existencebis} avec $ T =  \frac{\pi}{4} $. On applique à u la transformation de lentille définit en section 2.3 de \cite{moi} pour obtenir une fonction $ \tilde{u} $ qui, d'après les propositions 20 et 23 de \cite{moi}, vérifie les conditions du théorème. \hfill $ \boxtimes $
\begin{thm}  \label{2existence2}
Il existe $ s \in ]  \frac{d}{2}   - \frac{2}{p-1} ; \frac{d}{2} [ $, $ C_1,C_2 > 0 $ et $ \delta > 0 $ tels que pour tout $ 0 < T \leq 1 $, si $ u_0 \in E_0(\lambda ) $ avec $ \lambda < C_1 ( \arctan 2 T ) ^ {- \delta  }  $ alors il existe une solution à (\ref{2schrodinger}) sur $ [-T,T] $ dans l'espace $ e^{it\Delta}u_0 + B_{X_T^s } (0, C_2 \lambda ^p  ) $.
\end{thm}
\textit{Preuve.} Soit u donnée par le théorème \ref{2existencebis} à T remplacé par $ \frac{1}{2} \arctan 2T $. Puis, comme pour le théorème précédent, on applique à u la transformation de lentille définit en section 2.3 de \cite{moi} pour obtenir une fonction $ \tilde{u} $ qui, d'après la proposition 20 de \cite{moi} et la proposition \ref{2comparaisonbis}, vérifie les conditions du théorème. \hfill $ \boxtimes $
\\\\On démontre ensuite l'unicité des solutions construites.
\begin{thm} \label{2unicite}
Soient $ \frac{d}{2} > s >  \frac{d}{2}   - \frac{2}{p-1} $, $ u_0 \in E_0(\lambda) $ et $ T \in ]0,1 ] $. Supposons donné $ \tilde{u}_1 $ et $ \tilde{u}_2 $ deux solutions de (\ref{2schrodinger}) sur [-T,T] de l'espace $  e^{it\Delta} u_0 + X_T^s $ alors,
\begin{equation*}
 \tilde{u}_1(t) = \tilde{u}_2(t) \ \mbox{dans} \ L^2(\R^d), \ \forall t \in [-T,T].
\end{equation*} 
\end{thm}
\textit{Preuve.} Comme pour le théorème 69 de \cite{moi}, il suffit de prouver le théorème pour $ t \in [0,T] $.
\\Pour tout $ t \in \R $, on a 
\begin{align*}
& \partial_t || \tilde{u}_1(t) - \tilde{u}_2(t) ||^2_{ L^2(\R^d) } & 
\\ = \ & 2 \Re ( <  \partial_t(\tilde{u}_1(t)-\tilde{u}_2(t)) , \tilde{u}_1(t)-\tilde{u}_2(t)  >_{L^2(\R^d) \times L^2(\R^d) } )
\\ = \ & 2 |   <  | \tilde{u}_1(t)| ^{p-1} \tilde{u}_1(t) - |\tilde{u}_2(t)|^{p-1}  \tilde{u}_2(t)  , \tilde{u}_1(t)-\tilde{u}_2(t)  >_{L^2(\R^d) \times L^2(\R^d) } |
\\ \leq \ & 2 || \tilde{u}_1(t) - \tilde{u}_2(t) ||_{L^2(\R^d)} \times || \ |\tilde{u}_1(t)| ^{p-1} \tilde{u}_1(t) - |\tilde{u}_2(t)|^{p-1} \tilde{u}_2(t) ||_{L^2(\R^d)} 
\\ \leq \ & 2(p-1) ||\tilde{u}_1(t) - \tilde{u}_2(t) ||^2_{L^2(\R^d)} \times \left(  ||  \tilde{u}_1(t) ||^{p-1}_{L^\infty(\R^d)} +   ||  \tilde{u}_2(t) ||^{p-1}_{L^\infty(\R^d)}  \right).
\end{align*}
Puis, par le lemme de Gronwall, le théorème est prouvé si $ || \tilde{u}_1(t) ||^{p-1}_{L^\infty(\R^d)} +   ||  \tilde{u}_2(t) ||^{p-1}_{L^\infty(\R^d)} \in L^1_{loc} $ puisque $ || \tilde{u}_1(0) - \tilde{u}_2(0) ||^2_{ L^2(\R^d) } = 0 $. 
\\\\Mais, en utilisant les propositions \ref{2comparaisonbis} et \ref{2injection}, on obtient
\begin{align*}
|| \tilde{u}_i ||_{ L^{p-1}([0,T]), L^\infty(\R^d)  } & \leq || e^{it\Delta} u_0 ||_{ L^{p-1}([0,T]), L^\infty(\R^d) )  } +|| \tilde{v}_i ||_{ L^{p-1}([0,T]), L^\infty(\R^d) ) }
\\ & \leq C_T \times(     || e^{-itH} u_0 ||_{ L^{p-1}([-2 \pi,2 \pi ]), L^\infty(\R^d) )  } + || \tilde{v}_i ||_{ X_T^s } )
\\ & \leq C_T \times(    \lambda + ||  \tilde{v}_i ||_{ X_T^s } ).
\end{align*}
et le théorème est démontré. \hfill $ \boxtimes $
\\\\Enfin, on démontre que les solutions globales construites diffusent en $ \infty $ et en $ - \infty $.
\begin{thm} \label{2scattering}
Soit $ \tilde{u} $ l'unique solution globale de (\ref{2schrodinger}) construite dans le théorème \ref{2existence1} alors il existe $ L^+ \in \overline{H}^s(\R^d)  $ et $ L_- \in \overline{H}^s(\R^d) $ tels que
\begin{align*}
& \lim_{t \rightarrow \infty } || \tilde{u}(t) -e^{it\Delta} u_0 - e^{it\Delta} L^+ ||_{H^s(\R^d)} = 0 ,
\\ & \lim_{t \rightarrow - \infty } || \tilde{u}(t) -e^{it\Delta} u_0 - e^{it\Delta} L_- ||_{H^s(\R^d)} = 0. 
\end{align*}
\end{thm}
\textit{Preuve.} On pose $ T = \frac{\pi}{4}$, alors grâce aux propositions \ref{2estim1}, \ref{2estim2} et \ref{2estim3}, on obtient
\begin{align*}
 \int_0 ^ t  e ^{-i(t-s ) H }  K \cos (  2 s ) ^{  \frac{d(p-1)}{2} -2  }   & \left[ \ |   e^{-isH} u _ 0  + v(s) |^{p-1} \times (e^{-isH} u _ 0  + v(s) ) \right] ds 
 \\ & \in \overline{X}^s_T \hookrightarrow C^0([-T,T] , \overline{H}^s(\R^d) ).
\end{align*}
Ainsi, il existe $ L \in \overline{H}^s(\R^d) $ tel que  
\begin{equation*}
\lim_{t \rightarrow  T} \bigg| \bigg| \int_0 ^ t  e ^{-i(t-s) H } K \cos (  2 s ) ^{ \frac{d(p-1)}{2} -2  }   \left[ \ |   e^{-isH} u _ 0  + v(s) |^{p-1}  (e^{-isH} u _ 0  + v(s) ) \right] \ ds -  L \bigg| \bigg|_{ \overline{H}^s }  = 0,
\end{equation*}
puis
\begin{equation*}
\lim_{t \rightarrow  T} \bigg| \bigg| \int_0 ^ t  e ^{is H }  K \cos (  2 s ) ^{  \frac{d(p-1)}{2} -2  }   \left[ \ |   e^{-isH} u _ 0  + v(s) |^{p-1}  (e^{-isH} u _ 0  + v(s) ) \right] \ ds -  e^{iTH} L \bigg| \bigg|_{ \overline{H}^s  }  = 0.
\end{equation*}
Or, d'après le lemme 70 de \cite{moi}, on obtient
\begin{align*}
 \tilde{v} (t) & = \widetilde{  e^{-itH}  \int_0^t         e ^{is  H } K \cos (  2 s ) ^{  \frac{d(p-1)}{2} -2  }   \left[ \ |   e^{-isH} u _ 0  + v(s) |^{p-1} * (e^{-isH} u _ 0  + v(s) ) \right] ds             } 
 \\ & = e^{it\Delta} \int_0^{  \frac{1}{2} \arctan 2 t   } e ^{is  H } K \cos (  2 s ) ^{  \frac{d(p-1)}{2} -2  }   \left[ \ |   e^{-isH} u _ 0  + v(s) |^{p-1} * (e^{-isH} u _ 0  + v(s) ) \right] ds , 
\end{align*}
et donc
\begin{align*}
\lim _{t \rightarrow \infty } || \tilde{v} (t) - e^{it \Delta} e^{iTH} L  ||_{  H^s(\R^d) } = 0.
\end{align*}
\begin{flushright}
$ \boxtimes $
\end{flushright}
\section{Estimation de la régularité de la donnée initiale aléatoire}
\begin{dfn}
Pour $ t > 0 $, définissons
\begin{equation*}
\Omega_t = \left( \omega \in \Omega /  u_0 ^\omega  \in E_0( t )  \right).
\end{equation*}
\end{dfn}
Le but de cette partie est d'établir le théorème suivant:
\begin{thm} \label{2cond2} Sous les hypothèses (\ref{2hypothese1}) et (\ref{2hypothese2}) ou (\ref{2hypothese2bis}), il existe des constantes $ m ( \gamma ) , C, c > 0 $ telles que pour tout $ t > 0 $,
\begin{equation*} 
P(  \Omega^c_t  )\leq C \exp \left(  - c \left( \frac{  t}{   ||u_0||_{     \overline{H}^{ (d-1)/2  } (\R^d)  } } \right) ^{m(\gamma)} \right), 
\end{equation*}
où \begin{equation*}
m(\gamma) = \left\{
    \begin{array}{llll}
     &   \frac{2 \gamma}{2+\gamma}  & \mbox{ sous (\ref{2hypothese2}) } & \mbox{ si } \gamma \in ]0,1  ] ,   
    \\ &  \frac{3\gamma}{2\gamma+3}  & \mbox{ sous (\ref{2hypothese2bis}) } & \mbox{ si } \gamma  \in ]0,1  ] , 
    \\ &  \gamma    & \mbox{ sous  (\ref{2hypothese2bis})  } &  \mbox{ si } \gamma  \in ]1,2 ]          , 
  \\  &  2 & \mbox{ sous (\ref{2hypothese2bis}) }  &  \mbox{ si } \gamma  \geq 2  .
    \end{array}
\right.
\end{equation*}
\end{thm}
Par l'inégalité triangulaire, nous pouvons écrire
\begin{center}  $ P(  \Omega^c_t  ) $ \end{center}
\begin{equation} \label{2deuxtermes}
\leq P \left( \omega \in \Omega / \ ||u^\omega_0||_{     \overline{H}^{ (d-1)/2  }(\R^d) } \geq t \right)   +  P  \left( \omega \in \Omega / \ ||  e^{-itH} u_0 ||_{ L^{2p}( [-2\pi,2\pi] , \overline{W}^{  \frac{1}{7} , \infty } (\R^d))  }  \geq t  \right)
\end{equation}
et il suffit de montrer la majoration du théorème \ref{2cond2} pour chacun de ces deux termes. On commence par évaluer les moments de nos variables aléatoires à travers le lemme suivant :
\begin{lem} \label{2moments uniformes}
Sous l'hypothèse (\ref{2hypothese1}), il existe des constantes $  C_1,C_2,c > 0 $ telles que pour tout $ p \geq 1 $ et $ n \in \mathds{N} $,
\begin{align*}
 E( |g_n|^p) \leq  \left\{   \begin{array}{ll}
 C_1 \times  \left( \frac{p}{\gamma c } \right) ^{  \frac{p}{\gamma}  } & \mbox{ si } p \geq \gamma,
\\ C_2 & \mbox{ si } p \leq \gamma.
\end{array}
 \right.
\end{align*}
\end{lem}
\textit{Preuve.} On a
\begin{align*}
E( |g_n|^p) & = p \int_0^\infty \rho^{p-1} \times P \left( \omega \in \Omega / |g_n( \omega )| \geq \rho \right) \ d \rho
\\ & \leq p \int_0^\infty \rho^{p-1} \times Ce^{-c \rho ^\gamma} \ d \rho
\\ & \leq \frac{Cp}{\gamma}  \times \left( \frac{1}{c} \right) ^{ \frac{p}{\gamma} } \times  \int_0^\infty \mu^{  \frac{p-\gamma}{\gamma}} \times e^{-\mu} \ d \mu
\\ & \leq \frac{Cp}{\gamma}  \times \left( \frac{1}{c} \right) ^{ \frac{p}{\gamma} } \times  \Gamma \left(  \frac{p}{\gamma} \right) ,
\end{align*}
où $ \Gamma $ désigne la fonction gamma d'Euler. En utilisant les estimées de la fonction $ \Gamma $ suivantes :
\begin{align*} \begin{array}{ll}
 \Gamma(x) \leq c(Cx)^{x-1} & \mbox{ pour } x \geq 1,
\\  \Gamma(x) \leq \frac{C}{x} & \mbox{ pour } x \leq 1,
\end{array}
\end{align*}
on prouve le résultat. \hfill $ \boxtimes $
\\\\Puis, grâce à ce dernier lemme, nous pouvons estimer le premier terme de (\ref{2deuxtermes}).
\begin{prop} \label{2clef2} Sous l'hypothèse (\ref{2hypothese1}), il existe des constantes $ C,c > 0  $ telles que pour tout $ t >  0 $,
\begin{align*}
 P  \left( \omega \in \Omega \ / \   ||  u_0^\omega  ||_{   \overline{H}^{ (d-1)/2 } (\R^d)    }  \geq t \right)  \leq C \exp \left( -\frac{c t ^\gamma}{||u_0||^\gamma_{  \overline{H}^{ (d-1)/2 } (\R^d)    } } \right) .
\end{align*} 
\end{prop}
\textit{Preuve.} Il suffit d'établir l'estimation pour $  t \geq C ||u_0||_{  \overline{H}^{ (d-1)/2 } (\R^d)    }  $. Soit $ q \geq \max( 1 , \frac{\gamma}{2} ) $ alors d'après l'inégalité de Markov et le lemme \ref{2moments uniformes}, on trouve
\begin{align*}
P \left(  \omega \in \Omega \ / \   ||  u_0^\omega  ||_{   \overline{H}^{ (d-1)/2 } (\R^d)    }  \geq t \right) & = P   \left(  \omega \in \Omega \ / \  \sum_{n \in \mathds{N}} \lambda_n^{d-1} |c_n|^2 |g_n(\omega)|^2  \geq t^2 \right)
\\ & \leq   t^{-2 q} \times  E_{\omega} \left( \sum_{n \in \mathds{N}} \lambda_n^{d-1} |c_n|^2 |g_n(\omega)|^2   \right)  ^q 
\\ & \leq   t^{-2 q} \times  \bigg| \bigg| \sum_{n \in \mathds{N}} \lambda_n^{d-1} |c_n|^2 |g_n(.)|^2   \bigg| \bigg|_{L^q( \Omega)}^q 
\\ & \leq   t^{-2 q} \times  \left(  \sum_{n \in \mathds{N}} \lambda_n^{d-1} |c_n|^2 || g_n(.) ||^2_{L^{2q}( \Omega)}   \right)  ^q 
\\ & \leq     \left(  C \times  \frac{ ||u_0||_{ \overline{H}^{ \frac{d-1}{2} } ( \R^d ) }     }{t}  \times  \left( \frac{2q}{\gamma c } \right) ^{ \frac{1}{\gamma}  }  \right) ^{2q} .
\end{align*}
Puis, nous pouvons choisir $ q = \frac{\gamma c }{2} \times \left( \frac{ t }{2C||u_0||_{ \overline{H}^{ \frac{d-1}{2} } ( \R^d ) } } \right) ^\gamma \geq  \max( 1 , \frac{\gamma}{2} )  $ pour obtenir
\begin{align*}
P \left(  \omega \in \Omega \ / \   ||  u_0^\omega  ||_{   \overline{H}^{ (d-1)/2 } (\R^d)    }  \geq t \right) & \leq  \frac{1}{2^{2q}}
\\ & \leq  e^{ -2\ln(2)q  }  \leq \exp \left( -\frac{c t ^\gamma}{||u_0||^\gamma_{  \overline{H}^{ (d-1)/2 } (\R^d)    } } \right) .
\end{align*}
\begin{flushright}
$ \boxtimes $
\end{flushright}
Dès lors, il reste le second terme de (\ref{2deuxtermes}) à estimer. Pour cela, rappelons les estimées des fonctions propres de l'oscillateur harmonique dont la preuve peut être trouvée en corollaire 3.2 de \cite{koch}.
\begin{prop}  \label{2propre1}
Pour tout $ p \in [4 , \infty ]  $, il existe une constante $ C > 0 $ telle que pour tout $ n \in \mathds{N} $,
\begin{align*}
  &  || h_n ||_{L^p(\R^d)}  \leq C \lambda_n^{- \frac{1}{6} } & \mbox{ si }  & d =1,  &
 \\ &  || h_n ||_{L^p(\R^d)}  \leq C \lambda_n^{-1 + \frac{d}{2} } & \mbox{ si } &  d \geq 2.
\end{align*}
\end{prop}
On établit ensuite la proposition fondamentale suivante qui permet d'estimer le second terme de (\ref{2deuxtermes}).
\begin{prop} \label{2propositionfondamentale} On suppose qu'il existe une constante $ C> 0 $ telle pour toute suite $ (c_n)_{n \in \N} \in l^2(\N) $ et tout $ q \geq \max ( 2 , \gamma) $,
\begin{equation} \label{2hypothese5} \tag{$E_\gamma $}
\bigg| \bigg| \sum_{n \in \N} c_n g_n(\omega)  \bigg| \bigg|_{L^q(\Omega)} \leq C \times q ^{  \frac{1}{m(\gamma)} } \times \sqrt{ \sum_{n \in \N} |c_n|^2  } ,
\end{equation}
alors, sous cette condition, il existe deux constantes $ C,c > 0 $ telles que pour tout $ t \geq 0$,
\begin{equation*}
P \left(  \omega \in \Omega /  || e^{-itH} u_0^\omega  ||_{  L^{2p}( [-2 \pi, 2 \pi ]  ,  \overline{W}^{ \frac{1}{7}  , \infty  }(\R^d))   }  \geq t \right) \leq C  \exp \left( -  c \left( \frac{ t }{||u_0||_{  \overline{H}^{ (d-1)/2 } (\R^d)    } } \right) ^{  m(\gamma) } \right).
\end{equation*} 
\end{prop}
\textit{Preuve.} Comme 
\begin{equation*}
\overline{W}^{  \frac{1}{6} , r } ( \R^d )  \hookrightarrow \overline{W}^{  \frac{1}{7} , \infty  } (\R^d) \hspace*{1cm} \mbox{ si   }  \hspace*{1cm} \frac{1}{6} - \frac{1}{7} > \frac{d}{r},
\end{equation*}
on se ramène à démontrer l'existence de deux constantes $ C,c > 0 $ telles que pour tout $ t \geq 0 $ et $ r \geq 2 $,
\begin{equation}  \label{2terme2bis}
 P \left(  \omega \in \Omega /  || e^{-itH} u_0^\omega  ||_{  L^{2p}( [-2 \pi, 2 \pi ]  ,  \overline{W}^{ \frac{1}{6}  , r  }(\R^d))   }  \geq t \right) \leq C \exp \left( -  c \left( \frac{t}{||u_0||_{  \overline{H}^{ (d-1)/2 } (\R^d)    } } \right) ^{ m(\gamma) } \right) .
\end{equation} 
Il suffit de montrer l'estimation pour $ t \geq C ||u_0||_{ \overline{H}^\sigma(\R^d) } $. D'après les inégalités de Markov et Minkowsky, on obtient pour $ q \geq  \max( 2p,r , \gamma )$,
\begin{align*}
& P \left(  \omega \in \Omega \ / \  || e^{-itH} u_0^\omega  ||_{  L^{2p}( [-2 \pi, 2 \pi ]  ,  \overline{W}^{ \frac{1}{6}  , r  }(\R^d))   }  \geq t \right) 
\\ \leq \ & t^{-q } \times || e^{-itH} u_0^\omega  ||^q_{ L^q( \Omega, L^{2p}( [-2 \pi, 2 \pi ]  ,  \overline{W}^{ \frac{1}{6}  , r  }(\R^d)) )   }
\\  \leq \ & t^{-q } \times || H^{\frac{1}{12} } e^{-itH} u_0^\omega  ||^q_{ L^{2p}( [-2 \pi, 2 \pi ]  , L^ r  (\R^d , L^q( \Omega)) )   }.
\end{align*}
Puis, grâce à l'hypothèse (\ref{2hypothese5}), on obtient
\begin{align*}
|| H^{\frac{1}{12} } e^{-itH} u_0^\omega  ||^q_{ L^q( \Omega )   } & = \bigg| \bigg| \sum_{n \in \mathds{N}} \lambda_n^{\frac{1}{6} } c_n e^{-it\lambda_n^2} h_n(x) g_n(\omega )  \bigg| \bigg|^q_{ L^q( \Omega )   }
\\ & \leq C q ^{   \frac{1}{m(\gamma)} } \times  \sqrt{ \sum_{n \in \mathds{N}} \lambda_n^{\frac{1}{3} } |c_n|^2 .| h_n(x)|^2   }.
\end{align*}
Et finalement, par l'inégalité triangulaire et la proposition \ref{2propre1}, on a 
\begin{align*}
& P \left(  \omega \in \Omega \ / \  || e^{-itH} u_0^\omega  ||_{  L^{2p}( [-2 \pi, 2 \pi ]  ,  \overline{W}^{ \frac{1}{6}  , r  }(\R^d))   }  \geq t \right) 
\\ \leq \ & \left(  \frac{Cq^{ \frac{1}{m(\gamma)} }}{t}  \right)^q \times \bigg| \bigg| \sqrt{ \sum_{n \in \mathds{N}} \lambda_n^{\frac{1}{3} } |c_n|^2 .| h_n(x)|^2   } \bigg| \bigg|^q_{L^r(\R^d)}
\\ \leq \ & \left(  \frac{Cq^{ \frac{1}{m(\gamma)} }}{t}  \right)^q \times \bigg| \bigg|  \sum_{n \in \mathds{N}} \lambda_n^{\frac{1}{3} } |c_n|^2 .| h_n(x)|^2    \bigg| \bigg| ^{q/2}_{L^{r/2}(\R^d)}
\\ \leq \ & \left(  \frac{Cq^{ \frac{1}{m(\gamma)} }}{t}  \right)^q \times  \left( \sum_{n \in \mathds{N}} \lambda_n^{\frac{1}{3} } |c_n|^2 \ ||   h_n(x)   ||^2_{L^{r}(\R^d)} \right)^{q/2}
\\ \leq \ & \left(  \frac{Cq^{ \frac{1}{m(\gamma)} }  ||u_0||_{  \overline{H}^{(d-1)/2}(\R^d) } }{t}  \right)^q.
\end{align*}
Ainsi, il suffit de choisir $ q = \left( \frac{t}{2C||u_0||_{  \overline{H}^{(d-1)/2}(\R^d) } } \right)^{m(\gamma)} $ pour obtenir (\ref{2terme2bis}). \hfill $ \boxtimes$
\\\\Dès lors, on se ramène donc à démontrer (\ref{2hypothese5}) pour obtenir le théorème \ref{2cond2}.
\subsection{Preuve de (\ref{2hypothese5}) sous (\ref{2hypothese2}) si $ \gamma \in ]0,1  ] $}
Dans \cite{queffelec}, théorème 4.6, on a le lemme suivant :
\begin{lem} \label{2somme avec moments} Sous l'hypothèse (\ref{2hypothese2}), il existe une constante $ C> 0 $ telle que pour tout $ q \geq 1 $ et $ (c_n)_{n \in \mathds{N}} \in l^2(\N) $,
\begin{equation*}
E \left( \left(  \sum_{n \in \mathds{N}} c_n * g_n(\omega)  \right) ^{2q}  \right) \leq   (Cq)^q \times  E \left( \left(  \sum_{ n \in \mathds{N} }  \left( c_n * g_n(\omega) \right) ^2 \right)^q  \right) .
\end{equation*}
\end{lem}
Puis, en utilisant le lemme \ref{2moments uniformes}, on obtient
\begin{align*}
\bigg| \bigg|   \sum_{n \in \mathds{N}} c_n  g_n(\omega)    \bigg| \bigg|_{L^{2q}( \Omega  ) }^{2q} & \leq (Cq)^q \times \bigg| \bigg|   \sum_{n \in \mathds{N}} c_n^2  g_n(\omega)^2    \bigg| \bigg|^q_{L^{q}( \Omega  ) }
\\ & \leq (Cq)^q \times \left( \sum_{n \in \N} |c_n|^2 \times ||g_n||^2_{ L^{2q}(\Omega)  }    \right)^q
\\ & \leq \left( C \times q^{1+\frac{2}{\gamma} } \times \sum_{n \in \N} |c_n|^2   \right)^q.
\end{align*}
Et (\ref{2hypothese5}) est prouvé.
\subsection{Preuve de (\ref{2hypothese5}) sous (\ref{2hypothese2bis}) si $ \gamma \in ]0,1  ] $}
On commence par introduire une nouvelle définition :
\begin{dfn}
Pour $ p \in \mathds{N}^* $, on définit 
\begin{align*}
\mathcal{B}_{2p} = \bigg\{ \sigma \in \mathcal{S}_{2p} \ / \  \sigma(i) \neq i , \forall i \in \lbrace 1, .. , 2 p \rbrace \mbox{ et la longueur des} & \mbox{ cycles disjoints de } 
\\ & \sigma \mbox{ est égal à 2 ou 3 } \bigg\}.
\end{align*}
\end{dfn}
Puis, on établit quelques propriétés de $ \mathcal{B}_{2p}$.
\begin{lem} \label{2estimeebilineaireprobaesperancenulle1} Supposons donné $ (X_n)_{ n \in \mathds{N}} $ une suite de variables aléatoires vérifiant 
\\$ E(X_n)=0 $ et $ (n_1, ... n_{2p} )  \in \mathds{N}^{2p} $.
\\\begin{center}
Si $ E( X_{n_1} \times ... \times X_{n_{2p}}  ) \neq 0 $,
\end{center}
\begin{center}
alors, il existe $ \sigma \in \mathcal{B}_{2p} $ telle que $ n_{\sigma(i) } = n_i , \ \forall i \in \lbrace 1, .. , 2p\rbrace $.
\end{center}
\end{lem}
\textit{Preuve.} Ce résultat est clair par récurrence sur p. \hfill $ \boxtimes $
\begin{lem} \label{2estimeebilineaireprobaesperancenulle2}
Il existe une constante $ C> 0 $ telle que pour tout $ p \in \mathds{N}^* $, 
\\\begin{center} $ Card ( \mathcal{B}_{2p} ) \leq (Cp)^{ \frac{4p}{3} }$. \end{center}
\end{lem}
\textit{Preuve.} On utilise la formule de Stirling. On a 
\\-pour un ensemble à $ 2n $ éléments, le nombre de permutation qui ne fixe aucun point et constituée uniquement de transpositions dans leur décomposition en cycles disjoints est égal à $ \frac{(2n)!}{2^n n!} \leq (Cn)^{n}  $.
\\-pour un ensemble à $ 3n $ élément, le nombre de permutation qui ne fixe aucun point et constituée uniquement de 3-cycles dans leur décomposition en cycles disjoints est égal à $ \frac{(3n)!}{3^n n!} \leq (Cn)^{2n}  $.
\\\\Par conséquent, on obtient 
\begin{align*}
Card (\mathcal{B}_{2p} ) \leq C ^ p \times  \sum_{k=0}^p \binom{2p}{2k} \ k ^{ k  }  \ (p-k)^{\frac{4(p-k)}{3}} \leq (Cp)^{ \frac{4p}{3} }. & \hspace*{3cm} \boxtimes 
\end{align*}
Pour prouver (\ref{2hypothese5}), il suffit de traiter le cas $ q = 2 p $ où $ p \in  \mathds{N}^* $ et nous devons prouver que
\begin{equation*}
\bigg| \bigg|  \ \sum_{ n \in \mathds{N}   }  c_{n}  g_n(\omega)  \bigg| \bigg|^{2p}_{L^{2p}(\Omega)} \leq (Cp)^{ \frac{2p(2\gamma+3)}{3\gamma}  }   \times \left(  \sum_{n \in \mathds{N}} |c_{n} |^2 \right)^p.
\end{equation*}
On utilise les lemmes \ref{2estimeebilineaireprobaesperancenulle1} et \ref{2estimeebilineaireprobaesperancenulle2} pour obtenir que
\begin{align*}
\bigg| \bigg|  \ \sum_{ n \in \mathds{N}  }  c_{n}   g_n(\omega)  \bigg| \bigg|^{2p}_{L^{2p}(\Omega)} & =   \sum    \limits_{  { n_1,... ,n_{2p}  } } c_{n_1} ... c_{n_{p}} . \overline{ c_{n_{p+1}} }...\overline{c_{n_{2p}}} \times   E \left( \prod_{i=1}^{2p}  g_{n_i}   \right)
\\ & \leq   \sum    \limits_{  { n_1,...,n_{2p}  } }   | c_{n_1}| ...  |c_{n_{2p}}| \times \left|  E \left( \prod_{i=1}^{2p}  g_{n_i}   \right) \right|  
\\ & \leq   \bigg( \sum_{\sigma \in \mathcal{B}_{2p} }  \sum    \limits_{ \underset{n_{\sigma(i)} = n_i }{ n_1,...,n_{2p},  } }   | c_{n_1}| ...  |c_{n_{2p}}| \bigg) \times \  \underset{n \in \mathds{N}}{ \sup } \ E \left( |g_n|^{2p}  \right) 
\\ & \leq   Card (  \mathcal{B}_{2p} ) \times  (Cp)^{ \frac{2p}{\gamma} }   \times  \sup_{ \sigma \in \mathcal{B}_{2p}  } \bigg( \sum    \limits_{ \underset{n_{\sigma(i)} = n_i }{ n_1,...,n_{2p},  } }   | c_{n_1}| ...  |c_{n_{2p}}| \bigg)
\\ & \leq   (Cp)^{ \frac{2p(2\gamma+3)}{3\gamma}  }  \times   \left( \sum_{n \in \mathds{N}} |c_{n} |^2 \right)^p,
\end{align*}
où dans la dernière inégalité on utilise que $ l^2( \mathds{N} ) \hookrightarrow l^p( \mathds{N} )$ pour $ p \geq 2 $. \hfill $ \boxtimes $
\subsection{Preuve de (\ref{2hypothese5}) sous (\ref{2hypothese2bis}) si $ \gamma \in ]1,2] $}
Dans cette partie, on s'inspire de la preuve de \cite{burq4} en essayant de remplacer l'hypothèse ($H_2$) par l'hypothèse (\ref{2hypothese1}).
\begin{prop} \label{2majoration1} Soit $ X $ une variable aléatoire d'espérance nulle telle qu'il existe des constantes $ C,c> 0 $ telles que pour tout $ t \in [-2,2]$,
\begin{align*}
E(e^{|tX|}) \leq Ce^{ct^2},
\end{align*}
alors il existe une constante $ c> 0 $ telle que pour tout $ t \in [-1,1] $, 
\begin{equation*}
E(e^{tX}) \leq e^{ct^2}.
\end{equation*}
\end{prop}
\textit{Preuve.} De 
\begin{equation*}
e^u = 1+u+u^2 \int_0^1 (1-\theta ) e^{u \theta} \ d \theta,
\end{equation*}
on déduit pour tout $ t \in [-1;1]$ que
\begin{align*}
E( e^{tX}) & = 1 + t^2 \int_0^1 (1-\theta ) E( X^2 e^{t \theta X}  )  \ d \theta
\\ & \leq  1 + t^2 \int_0^1  \sqrt{  E( X^4 ) E(  e^{2t \theta X}  )  } \ d \theta
\\ & \leq  1 + C t^2 \int_0^1     e^{2c t^2 \theta^2 }     \ d \theta
\\ & \leq  1 + C t^2 e^{2c t^2 }  
\\ & \leq  e^{c' t^2 }.
\end{align*}
\begin{flushright}
$ \boxtimes $
\end{flushright}
\begin{prop} \label{2majoration2} Sous les hypothèses (\ref{2hypothese1}) et (\ref{2hypothese2bis}), il existe une constante $ c> 0 $ telle que pour tout $ n \in \mathds{N} $ et $ t \in \R $,
\begin{align*}
E(e^{tg_n} ) \leq \left\{
    \begin{array}{ll}
        e^{ct^2} & \mbox{ si } |t| \leq 1, \\
        e^{c|t|^{  \frac{\gamma}{\gamma-1} }} & \mbox{ si } |t|\geq 1.
    \end{array}
\right.
\end{align*}
\end{prop}
\textit{Preuve.} Sous l'hypothèse (\ref{2hypothese1}), on a pour tout $ t \in \R$,
\begin{align*}
E(e^{|t||g_n|} ) & = 1 + |t| \times \int_0^\infty e^{|t| \rho} \times P \left( \omega \in \Omega / |g_n(\omega )| \geq \rho  \right) \ d \rho 
\\ & \leq  1 + C |t| \times \int_0^\infty e^{  |t| \rho-c|\rho|^\gamma}  \ d \rho 
\\ & \leq  1 + C |t| \times \sup_{ \rho \in \R^+} \left( e^{  |t| \rho- \frac{c}{2} |\rho|^\gamma} \right) \times  \int_0^\infty e^{  -\frac{c}{2} |\rho|^\gamma}  \ d \rho 
\\ & \leq  1 + C(\gamma) \times  |t| \times  e^{  |t|^{  \frac{\gamma}{\gamma-1} } \times \left( \frac{2}{c \gamma} \right)^{1/(\gamma-1)} \times \frac{\gamma-1}{\gamma} } 
\\ & \leq  1 + C(\gamma) \times  |t| \times  e^{ c(\gamma) |t|^{  \frac{\gamma}{\gamma-1} }  }.
\end{align*}
Par conséquent, 
\begin{align*}
E(e^{|t||g_n|} ) \leq \left\{
    \begin{array}{ll}
        C' e^{c't^2} & \mbox{ si } |t| \leq 2, \\
        e^{c'|t|^{  \frac{\gamma}{\gamma-1} }} & \mbox{ si } |t|\geq 1,
    \end{array}
\right.
\end{align*}
et nous pouvons utiliser la proposition \ref{2majoration1} pour conclure. \hfill $ \boxtimes $
\begin{prop} \label{2majoration3} Sous les hypothèses (\ref{2hypothese1}) et (\ref{2hypothese2bis}), il existe deux constantes $ C,c> 0 $ telles que pour tout $ \rho \geq 0 $ et $ (c_n)_{n \in \mathds{N} } \in l^2(  \N) $,
\begin{equation*}
P \left(  \omega  \in \Omega / \bigg|  \sum_{n \in \mathds{N}} c_n g_n(\omega ) \bigg| \geq \rho \right) \leq C e^{ -c \left( \frac{\rho}{ ||c_n||_{l^2(\N)}  }    \right)^\gamma }.
\end{equation*}
\end{prop}
\textit{Preuve.} On écrit,
\begin{align*}
& P \left(  \omega  \in \Omega / \bigg|  \sum_{n \in \mathds{N}} c_n g_n(\omega ) \bigg| \geq \rho \right) \\ \leq \ & P  \left(  \omega  \in \Omega /   \sum_{n \in \mathds{N}} c_n g_n(\omega )  \geq \rho \right) + P \left(  \omega  \in \Omega /   \sum_{n \in \mathds{N}} -c_n g_n(\omega )  \geq \rho \right),
\end{align*}
et il suffit de montrer la majoration pour le premier terme.
\\\\D'après l'inégalité de Markov et la proposition \ref{2majoration2}, on obtient pour tout $ t \geq 0 $,
\begin{align*}
P \left(  \omega  \in \Omega /   \sum_{n \in \mathds{N} } c_n g_n(\omega )  \geq \rho \right) & = P \left(  \omega  \in \Omega / \exp \left(  t \times  \sum_{n \in \mathds{N}}  c_n g_n(\omega )   \right) \geq e^{t\rho} \right)
\\ & \leq e^{-t\rho} \times E \left( \exp \left(  t \times  \sum_{n \in \mathds{N}}  c_n g_n(\omega ) \right) \right)
\\ & \leq e^{-t\rho} \times \prod_{n \in \mathds{N} } E \left( e^{   t \times  c_n g_n(\omega ) } \right)
\\ & \leq e^{-t\rho} \times \prod_{n \in \mathds{N} } \max \left( e^{c \times  |c_n . t |^2 }  , e^{c \times  |t .c_n|^{ \frac{\gamma}{\gamma-1}  } } \right)
\\ & \leq e^{-t\rho} \times  \exp \left( c \times \sum_{n \in \mathds{N} } |c_n . t |^2 \right)  \times  \exp \left( c \times \sum_{n \in \mathds{N} } |t .c_n|^{ \frac{\gamma}{\gamma-1}  } \right) 
\\ & \leq e^{-t\rho} \times  \exp \left( c \times ( t ||c_n||_{l^2} ) ^2 \right)  \times  \exp \left( c \times ( t ||c_n||_{l^2} )^{ \frac{\gamma}{\gamma-1}  } \right),
\end{align*}
puis, il suffit de choisir $ t = \epsilon \frac{\rho^{\gamma-1}}{  ||c_n||^\gamma_{l^2}  } $ pour obtenir le résultat souhaité. \hfill $ \boxtimes $
\begin{prop} \label{2majoration4} Sous les hypothèses (\ref{2hypothese1}) et (\ref{2hypothese2bis}), il existe une constante $ C> 0 $ telle que pour tout $ q \geq 2 $ et $ (c_n)_{n \in \mathds{N} } \in l^2(\N)$,
\begin{equation*}
\bigg| \bigg|  \sum_{n \in \mathds{N}} c_n g_n(\omega ) \bigg| \bigg|_{L^q(\Omega)} \leq C \times  q^{  \frac{1}{\gamma} }  \times \sqrt{  \sum_{n \in \N} |c_n|^2  } .
\end{equation*}
\end{prop}
\textit{Preuve.} En utilisant la proposition \ref{2majoration3}, on obtient
\begin{align*}
\bigg| \bigg|  \sum_{n \in \mathds{N}} c_n g_n(\omega ) \bigg| \bigg|^q_{L^q(\Omega)} & = q \int_0^\infty \rho^{q-1} \times P \left( \omega \in \Omega / \bigg|  \sum_{n \in \mathds{N}} c_n g_n(\omega ) \bigg| \geq \rho \right) \ d \rho
\\ & \leq  Cq \int_0^\infty \rho^{q-1} \times e^{ -c \left( \frac{\rho}{ ||c_n||_{l^2(\N)}  }    \right)^\gamma } \ d \rho
\\ & \leq  (C ||c_n||_{l^2(\N)} ) ^q \times q \int_0^\infty u^{q/\gamma-1} \times e^{ - u  }  \ d u
\\ & \leq  (C' ||c_n||_{l^2(\N)} ) ^q \times q^{ \frac{q}{\gamma}  }.
\end{align*}
Ce qui prouve la proposition. \hfill $ \boxtimes $
\\Et démontre (\ref{2hypothese5}).
\section{Preuves des théorèmes}
Dans cette section, on démontre les théorèmes \ref{2thm1}, \ref{2thm3} et \ref{2thm2}.
\subsection{Preuve du théorème \ref{2thm1}}
En utilisant les théorèmes \ref{2existence1}, \ref{2unicite}, \ref{2scattering}, pour obtenir le théorème \ref{2thm1}, il suffit d'établir que pour tout $ t > 0 $,
\begin{equation} \label{2cond1}
P( \Omega_t  ) > 0.
\end{equation}
Pour cela, introduisons la définition suivante :
\begin{dfn}
Pour $ u_0 =   \displaystyle{  \sum_{n \in \mathds{N}} c_n h_n(x) }  $ une fonction de $ L^2 ( \R^d ) $, nous définissons pour $ N \in \N^*$,
\begin{align*}
& [u_0]_N = \sum_{\lambda_n < N } c_n h_n(x),
\\ & [u_0]^N = \sum_{\lambda_n   \geq N} c_n h_n(x).
\end{align*}
\end{dfn}
\begin{prop} Sous les hypothèses (\ref{2hypothese1}) et (\ref{2hypothese2}) ou (\ref{2hypothese1}) et (\ref{2hypothese2bis}), pour tout $ t > 0 $, il existe $ N \in \N^ * $ tel que
\begin{align*}
& P(  \Omega _t )  \geq \frac{1}{2}  
\\ & \times P \left( \ \omega \in \Omega \ /  \  || \ [u^\omega_0]_N ||_{  \overline{H}^{(d-1)/2}  (\R^d) }  \leq \frac{t}{2}  \cap ||  e^{-itH} [u^\omega_0]_N ||_{ L^{2p}( [-2\pi,2\pi] , \overline{W}^{  \frac{1}{7} , \infty } (\R^d))  \ }  \leq \frac{t}{2} \right) . 
\end{align*}
\end{prop}
\textit{Preuve.} Par indépendante, en utilisant le théorème \ref{2cond2}, on obtient
\begin{align*}
&  P(  \Omega _t  )  
\\  = \ & P \bigg( \ \omega \in \Omega \ /  \  || \ [u^\omega_0]_N + [u^\omega_0]^N ||_{  \overline{H}^{(d-1)/2}  (\R^d) }  \leq t 
\\ &  \hspace*{2cm} \cap  || e^{-itH} [u^\omega_0]_N + e^{-itH} [u^\omega_0]^N  ||_{  L^{2p}( [-2\pi,2\pi] , \overline{W}^{ \frac{1}{7} , \infty  } (\R^d)) } \leq  t \bigg) 
\\ \geq \ & P \left( \ \omega \in \Omega \ /  \  || \ [u^\omega_0]_N ||_{   \overline{H}^{(d-1)/2} (\R^d)  } \leq \frac{t}{2}     \cap  || e^{-itH} [u^\omega_0]_N ||_{  L^{2p}( [-2\pi,2\pi] , \overline{W}^{ \frac{1}{7} , \infty  } (\R^d)) } \leq  \frac{t}{2} \right)  
\\ &  \hspace*{0.15cm} \times P \left( \ \omega \in \Omega \ /  \  || \ [u^\omega_0]^N ||_{  \overline{H}^{(d-1)/2}   (  \R^d ) } \leq \frac{t}{2}     \cap  || e^{-itH} [u^\omega_0]^N ||_{  L^{2p}( [-2\pi,2\pi] , \overline{W}^{ \frac{1}{7} , \infty  } (\R^d)) } \leq  \frac{t}{2} \right)  
\\ \geq \ & P \left( \ \omega \in \Omega \ /  \  || \ [u^\omega_0]_N  ||_{  \overline{H}^{(d-1)/2}   (  \R^d ) } \leq \frac{t}{2}     \cap  || e^{-itH} [u^\omega_0]_N ||_{  L^{2p}( [-2\pi,2\pi] , \overline{W}^{ \frac{1}{7} , \infty  } (\R^d)) } \leq  \frac{t}{2}  \right)
\\ & \hspace*{5.2cm}  \times \left( 1-  C \exp \left(  - \frac{c t ^{m(\gamma)}}{   (  \underset{\lambda_n \geq N  }{\sum}   \lambda_n ^{d-1} \times |c_n |^2  )^{m(\gamma)/2} }   \right)  \right) .
\end{align*}
Or
\begin{equation*}
\lim_{N \longrightarrow \infty } 1-  C \exp \left(  - \frac{c t ^{m(\gamma)}}{   (  \underset{\lambda_n \geq N  }{\sum}   \lambda_n ^{d-1} \times |c_n |^2  )^{m(\gamma)/2} }   \right) = 1,
\end{equation*}
ainsi, il existe $ N \in \N^* $ tel que $ 1-  C \exp \left(  - \frac{c t ^{m(\gamma)}}{   (  \underset{\lambda_n \geq N  }{\sum}   \lambda_n ^{d-1} \times |c_n |^2  )^{m(\gamma)/2} }   \right) \geq \frac{1}{2} . $ \hfill $ \boxtimes $
\\\\ Par conséquent, pour prouver (\ref{2cond1}), il suffit de prouver la proposition suivante :
\begin{prop} \label{2cond3} Sous l'hypothèse (\ref{2hypothese3}), pour tout $ t > 0 $ et $ N \in \N ^ * $,
\begin{align*} 
P \left( \ \omega \in \Omega \ /  \  || \ [u^\omega_0]_N ||_{  \overline{H}^{ (d-1)/2 }  (\R^d) }  \leq t    \cap  || e^{-itH} [u^\omega_0]_N  ||_{  L^{2p}( [-2\pi,2\pi] , \overline{W}^{ \frac{1}{7} , \infty  } (\R^d)) } \leq  t \right) > 0 .
\end{align*}
\end{prop}
\textit{Preuve.} En utilisant l'hypothèse (\ref{2hypothese3}), on obtient
\begin{align*}
& P \left( \ \omega \in \Omega \ /  \  || \ [u^\omega_0]_N ||_{  \overline{H}^{ (d-1)/2 }  (\R^d) }  \leq t    \cap  || e^{-itH} [u^\omega_0]_N  ||_{  L^{2p}( [-2\pi,2\pi] , \overline{W}^{ \frac{1}{7} , \infty  } (\R^d)) } \leq  t \right)
\\ \geq \ & P \left( \ \omega \in \Omega \ /  \   \sum_{ \lambda_n  < N   } \lambda_n^{d-1} |c_n|^2  | g_n(  \omega  )|^2  \leq   \frac{ C t ^2 }{N^{2d}} \right)
\\  \geq \ & P \left( \ \omega \in \Omega \ /  \   \sum_{ \lambda_n < N   }  | g_n(  \omega  )|^2  \leq   \frac{ C t^2  }{ N^{2d} \times || u _0 ||^2_{ \overline{H}^{d-1}(\R^d)   }  } \right)
\\  \geq \ & P   \left(  \underset{\lambda_n < N }{\bigcap}  \left( \ \omega \in \Omega \ /  \  | g_n(  \omega  )|^2  \leq   \frac{ C  t^2 }{N^{4d} \times || u _0 ||^2_{ \overline{H}^{d-1}(\R^d)   }  } \right) \right)
\\  \geq \ & \prod_{\lambda_n < N} P  \left( \ \omega \in \Omega \ /  \  | g_n(  \omega  )|^2  \leq \frac{ C  t^2 }{N^{4d} \times || u _0 ||^2_{ \overline{H}^{d-1}(\R^d)   }  } \right)  > 0. & \boxtimes 
\end{align*}
\subsection{Preuve du théorème \ref{2thm3}}
On adapte ici la preuve du paragraphe 5 de \cite{burq4}. Grâce aux théorèmes \ref{2existence2} et \ref{2unicite}, on sait que si $ u_0 \in E_0( (\arctan 2T) ^{-\delta} ) $ alors il existe une unique solution à l'équation (\ref{2schrodinger}) sur $ [-T,T] $ dans l'espace $ e^{it\Delta} u_0 + B_{ X^s_T } ( 0, C_T )  $.
\\\\Définissons
\begin{equation*}
\Omega_T = \left( \omega \in \Omega / u_0^\omega  \in  E_0( (\arctan 2T) ^{-\delta} )  \right),
\end{equation*}
alors par le théorème \ref{2cond2},
\begin{equation*}
P ( \Omega_T ^c ) \leq C \exp(   - c  ( \arctan 2T) ^{-\delta'}  )  ) .
\end{equation*}
Par conséquent, si nous posons
\begin{equation*}
\Sigma =  \underset{n \in \N^*}{ \cup }  \Omega_{ 1/n } 
\end{equation*}
alors $ P( \Sigma  ) = 1 $ et le théorème \ref{2thm3} est prouvé. \hfill $ \boxtimes $
\subsection{Preuve du théorème \ref{2thm2}}
Grâce aux théorèmes \ref{2existence1}, \ref{2unicite} et \ref{2scattering}, pour prouver le théorème \ref{2thm2}, il suffit d'établir que pour tout $ \lambda > 0 $,
\begin{equation} \label{2cond4}
\underset{\eta \rightarrow 0 }{ \lim } \ P  \left( \omega \in \Omega_\lambda^c | \ ||u_0^\omega||_{ \overline{H}^{(d-1)/2} ( \R^d )  } \leq \eta \right) = 0 .
\end{equation}
Nous pouvons ensuite utiliser la même méthode que la proposition A.7 de \cite{burq7} pour obtenir que 
\begin{equation*}
  P  \left( \omega \in \Omega^c_\lambda   \ | \ ||u_0^\omega||_{ \overline{H}^{(d-1)/2}(\R^d) }  \leq \eta \right)  \leq  C e ^{  -c \frac{\lambda^2}{\eta^2} }, 
\end{equation*}
et (\ref{2cond4}) est démontré.
\newpage
\nocite{*}
\bibliographystyle{short}
\bibliography{biblioarticle2}
\end{document}